\documentclass[10pt]{article}
\usepackage[english]{babel}
\usepackage[latin1]{inputenc} 
\usepackage{amssymb,amsmath}
\usepackage{hyperref}

\newtheorem{teor}{Theorem}[section]
\newtheorem{lema}[teor]{Lemma}
\newtheorem{ej}[teor]{Example}
\newtheorem{prop}[teor]{Proposition}
\newtheorem{corol}[teor]{Corollary}
\newtheorem{defin}[teor]{Definition}
\newtheorem{obs}[teor]{Remark}

\title{Distinguished connections on $(J^{2}=\pm 1)$-metric manifolds}
\author{Fernando Etayo\footnote{Departamento de Matem\'{a}ticas, Estad\'{\i}stica y Computaci\'{o}n. Facultad de Ciencias.  Universidad de Cantabria. Avda. de los Castros, s/n, 39071 Santander, SPAIN. e-mail: etayof@unican.es}\, and Rafael Santamar\'{\i}a\footnote{Departamento de Matem\'{a}ticas. Escuela de Ingenier\'{\i}as Industrial e Inform\'{a}tica. Universidad de Le\'{o}n. Campus de Vegazana, 24071 Le\'{o}n, SPAIN. e-mail: rsans@unileon.es}}

\date{\today}

\begin{document}

\maketitle

\begin{abstract}
We study several linear connections (the first canonical, the Chern, the well adapted, the Levi Civita, the Kobayashi-Nomizu, the Yano, the Bismut and those with totally skew-symmetric torsion)  which can be defined on the four geometric types of $(J^2=\pm1)$-metric manifolds. We characterize when such a connection is adapted to the structure, and obtain a lot of results about coincidence among connections. We prove that the first canonical and the well adapted connections define a one-parameter family of adapted connections, named canonical connections, thus extending to almost Norden and almost product Riemannian manifolds the families introduced in almost Hermitian and almost para-Hermitian manifolds in \cite{gauduchon} and \cite{ivanov}. We also prove that every  connection studied in this paper is a canonical connection, when it  exists and it is an adapted connection.
\end{abstract}



\section{Introduction}

In the present paper we study connections defined on manifolds having an $(\alpha,\varepsilon)$-structure. A manifold will be called to have an $(\alpha, \varepsilon)$-structure if $J$ is an almost complex ($\alpha =-1$) or almost product ($\alpha =1$) structure and $J$ is an isometry ($\varepsilon =1$) or anti-isometry ($\varepsilon =-1$). It is also said that $(M,J,g)$ is a $(J^{2}=\pm 1)$-metric manifold. Thus, there exist four kinds of $(\alpha ,\varepsilon )$ structures according to the values $\alpha ,\varepsilon \in \{-1,1\}$, where
\[
J^2 = \alpha Id, \quad g(JX,JY)= \varepsilon g(X,Y), \quad \forall X, Y \in \mathfrak{X}(M). 
\]
As is well known, these four geometries have been intensively studied. The corresponding manifolds are known as:

\begin{enumerate}
\renewcommand*{\theenumi}{\roman{enumi})}
\renewcommand*{\labelenumi}{\theenumi}

\item  Almost Hermitian manifold if it has a  $(-1,1)$-structure. We shall consider through this paper the case $g$ being a Riemannian metric.

\item Almost anti-Hermitian or almost Norden manifolds if it has a $(-1,-1)$-structure. The metric $g$ is semi-Riemannian having signature $(n,n)$.

\item Almost product Riemannian manifolds if it has an $(1,1)$-structure. We shall consider through this paper the case $g$ being a Riemannian metric and the trace of $J$ vanishing, which in particular means these manifolds have even dimension.

\item Almost para-Hermitian manifolds if it has an $(1,-1)$-structure. The metric $g$ is semi-Riemannian having signature $(n,n)$.
\end{enumerate}

\noindent If the structure $J$ is integrable, i.e., the Nijenhuis tensor $N_{J}=0$, the corresponding manifolds are called Hermitian, Norden, Product Riemannian and para-Hermitian (without the word ``almost"). Integrabilty means $M$ is a holomorphic manifold in cases $i)$ and $ii)$, and $M$ has two complementary foliations in cases $iii)$ and $iv)$.

A linear connection is said to be reducible, natural or adapted to a manifold $M$ having an $(\alpha,\varepsilon)$-structure $(J,g)$ if its covariant derivative $\nabla^{\mathrm{a}}$ parallelizes both structures, i.e., $\nabla^{\mathrm{a}} J=0, \nabla^{\mathrm{a}} g=0$. The most significative natural connection is the well adapted connection $\nabla^{\mathrm{w}}$, which has been intensively studied in \cite{racsam}. We say that it is the most significative connection because it measures the integrability of the $G_{(\alpha,\varepsilon )}$-structure defined by $(J,g)$: it is integrable if and only if the torsion and the curvature tensors of the well adapted connection vanish. Nevertheless there exist others connections on $M$ having very interesting properties, although they are not adapted in the general case. The star is the Levi Civita connection: it is uniquely defined in such a manifold (definition depends just on the metric) but it is adapted to the $(\alpha,\varepsilon)$-structure $(J,g)$ if and only if the manifold is of K\"{a}hler type. Many of the results obtained in the four geometries are expressed in terms of the Levi Civita connection.

We can say that both connections, the well adapted and the Levi Civita ones, are distinguished connections. There exist a plethora of connections which have been defined in some of these manifolds.  In the celebrated paper of Gauduchon \cite{gauduchon}, he wrote in the introduction: ``I propose to the Reader as a kind of \textit{vade mecum} for some basics of almost Hermitian geometry", including ``a unified presentation of a canonical class of (almost) Hermitian connections". Our purpose is to extend that unified presentation to all the four geometries,  describing different connections appearing in the Literature, obtaining relations among them, and extending results from some of the four geometries to the rest of them. Thus, there is a sensible amount of new results in the present paper, which will be showed later.  The following ideas are important through the paper: 
\begin{enumerate}
\item The most important number  is the product $\alpha \varepsilon$, thus existing two classes of structures: the first one ($\alpha \varepsilon =-1$) defined by almost Hermitian and almost para-Hermitian  and the second one ($\alpha \varepsilon =1$) by almost Norden and almost product Riemannian. As a basic example one can consider  the fundamental tensor $\Phi (X,Y)$ defined as $\Phi (X,Y)=g(JX,Y)$, $\forall X, Y \in \mathfrak{X} (M)$: it is a 2 form if $\alpha \varepsilon =-1$ and a metric if $\alpha \varepsilon =1$.
\item Some connections obtained in the Literature are given with an explicit formula while others are given by imposing a condition they satisfy (in general about the torsion tensor). In this second case, one should prove existence and uniqueness of such a connection.
\item Some of the connections are natural or adapted to the $G_{(\alpha ,\varepsilon)}$-structure while some other important connections are not adapted.
\item The more geometric properties the manifold has, the less number of different distinguished connections exist (in the case of K\"{a}hler type manifolds all distinguished connections studied here coincide).
\end{enumerate}

\noindent  These  ideas will be carefully showed through the paper. The paper is as self-contained as possible. Some of the known results are  proved again according to the approach to the topic given in the present paper. We do not follow a chronological order, this is not a historical or survey paper, but a logical order from less to more structure. The almost Hermitian geometry is the model for the other three geometries, but we show the results for all the four geometries together, as possible, thus offering a new perspective from which to view all the structures.  Many mathematicans have studied connections on a specific kind of $(J^{2}=\pm 1)$-metric manifolds. Main contributions we have followed are in papers quoted in References. We  classified them according to the geometry they study:
\begin{itemize}
\item Almost Hermitian: Bismut \cite{bismut};  Davidov, Grantcharov and Mu\u{s}karov \cite{davidov}; Ganchev and Kassabov \cite{ganchev-kassabov}; Gauduchon \cite{gauduchon}; Gray \cite{gray}; Gray and Hervella \cite{gray-hervella}; Rod Gover and Nurowski \cite{rod-nurowski}; Vezzoni \cite{vezzoni}. 
\item Almost para-Hermitian: Chursin, Sch\"{a}fer and Smoczyk \cite{chursinetal}; Cruceanu and   Etayo \cite{etayo};   Gadea and   Mu\~{n}oz Masqu\'{e} \cite{gadea};  Ivanov and Zamkovoy \cite{ivanov}; Olszak \cite{olszak}.
\item Almost Norden: Ganchev and  Borisov \cite{ganchev-borisov}; Ganchev and Mihova \cite{ganchev-mihova}; Mekerov \cite{mekerov2};  Mekerov and Manev \cite{mekerov-manev};   Teofilova \cite{teofilova2} and \cite{teofilova}.
\item Almost Product Riemannian with vanishing trace: Gribacheva and Mekerov \cite{gribacheva}; Mekerov \cite{mekerov}; Mihova \cite{mihova}; Staikova and Gribachev \cite{staikova}; Yano \cite{yano}.
\end{itemize}

\noindent Other references we have used can not be included in this elemental scheme. So Agricola \cite{agricola} and Friedrich and Ivanov \cite{friedrich} pay particular attention to non-integrable $G$-structures on Riemannian manifolds, where connections with totally skew-symmetric torsion, if there exist, play an important role in the study this kind of $G$-structures.  A classical and seminal reference about this topic is the book of  Lichnerowicz \cite{lichnerowicz2}. Some similar comments can be said about the book of Kobayashi and Nomizu \cite{KN}, where  basic results about almost complex and almost Hermitian geometry are stated.

In \cite{debrecen} we have studied in a unified way the geometric properties of  $(J^{2}=\pm 1)$-metric manifolds. In the recent paper \cite{racsam} we introduce the well adapted connection of any  $(J^{2}=\pm 1)$-metric manifold, thus being  our first approach to this unified vision of connections in the four geometries.
\bigskip

The organization of the paper is as follows:

Section 2 is devoted to the study of reducible connections on $(J^{2}=\pm 1)$-manifolds, i.e., manifolds having an almost complex or paracomplex structure, without a metric. We say they have an $\alpha$-structure. We are inspired in the works \cite{etayo} and \cite{rod-nurowski} which take this starting point in their study of para-Hermitian and Hermitian geometries, respectively. We obtain a characterization of  reducible connections (Lemma \ref{teor:JconexionesS} and Proposition \ref{teor:JconexionesQ}). A key point is the definition (see formulas (\ref{eq:Jconexion0}) and (\ref{eq:Jconexion1})) of two adapted  connections, whose covariant derivatives are denoted as $\nabla^{0}$ and $\nabla^{1}$. These connections depend on the selection of an arbitrary connection $\nabla$. They allow to parametrize all the natural connections (Proposition \ref{teor:JconexionesQ}) and the line they define $\{ (1-s)\nabla^{0} +s\nabla^{1} \colon s \in \mathbb{R}\}$ is formed by natural connections (Proposition \ref{teor:uniparametrica-alpha}). Besides we present Kobayashi-Nomizu and Yano type connections. All the connections in this Section depend on the selection of a connection $\nabla$ which is the basis of the definition of the other ones. 

In Section 3 we show a first collection of results about reducible connections on $(J^{2}=\pm 1)$-metric manifolds, i.e., manifolds $(M,J,g)$ having an $(\alpha ,\varepsilon )$-structure. In the previous paper (\cite[Lemma 4.3]{racsam}) we have parametrized the set of natural connections of such a structure, by means of the Levi Civita connection and the potential tensor (which is the difference tensor between a natural connection and the Levi Civita connection). In the present paper we introduce the first canonical connection $\nabla ^{0}$ (Definition \ref{teor:first-connection}) from the Levi Civita connection $\nabla ^{\mathrm{g}}$, following the ideas of the above Section. Thus, the connection $\nabla ^{0}$ is uniquely determined, and it is always natural (Lemma \ref{teor:nabla0natural}). In Lemma \ref{teor:naturalQ} we parametrize the set of natural connections taking $\nabla ^{0}$ as starting point.

The following Section 4 has technical character. We study some tensors derived from $\nabla^{\mathrm{g}} J$ which will be useful in the study of connections in the remaining sections. Given a $(J^2=\pm1)$-metric manifold $(M,J,g)$, we  study three tensors: the covariant derivative $\nabla^{\mathrm{g}} \Phi$ of  the fundamental tensor $\Phi$, the Nijenhuis tensor $N_{J}$ and the second Nijenhuis tensor $\widetilde N_J^{\alpha\varepsilon}$. Properties of the two first tensors are well known and will be summarized. We focus on the expression of all of these tensors by means of $\nabla^{\mathrm{g}} J$. The relations we obtain between Nijenhuis and torsion tensors allow us to obtain sufficient conditions for the integrability of $J$ expressed by means of the torsion of an adapted connection to $(M,J,g)$ (Lemmas \ref{teor:nijenhuis-torsion} and \ref{teor:mascondicionessuficientes}). The vanishing of the second Nijenhuis tensor characterizes quasi-K\"{a}hler manifolds, as we prove in Propositions \ref{teor:caracterizacion-quasikahler+1} and \ref{teor:caracterizacion-quasikahler-1}. The expression of this tensor depends on the value $\alpha \varepsilon =-1$ or $\alpha \varepsilon =1$ (see  Definition \ref{teor:njae} and formulas (\ref{eq:njt}) and (\ref{eq:njd})). For this reason we need two different characterization theorems. In the  case $\alpha \varepsilon =1$ this tensor was known (see, {\em e.g.}; \cite{ganchev-borisov} and \cite{staikova}) but as far as we know there was no a definition for $\alpha \varepsilon =-1$. Last results in this Section provide relations among the vanishing of the quoted tensors and the type of the manifold.

Sections 5 and 6 are the core of the paper. In Section 5 we study the following distinguished connections on a $(J^{2}=\pm 1)$-metric manifold $(M,J,g)$: the first canonical, the Chern (which can be defined just in the case $\alpha\varepsilon=-1$), the well adapted, the Levi Civita, the Kobayashi-Nomizu, the Yano and those with totally skew-symmetric torsion. The three first  are always natural connections. In the case of the remaining four, one needs to determine the conditions to be satisfied in order to be natural (which are summarized in Table \ref{table:connections}). The principal difference which allows to group them is the following: the first canonical, the Kobayashi-Nomizu and the Yano connections are uniquely defined from the Levi Civita connection whilst the Chern, the well adapted connection and those with totally skew-symmetric torsion are defined imposing a condition about the torsion. We start the Section relating the torsion tensor of the first canonical connection with the integrability of the $\alpha$-structure $J$ (Corollary \ref{teor:toro-integrability}). We follow recalling the unified presentation of the Chern connection in the case $\alpha\varepsilon=-1$ obtained in \cite{racsam} (Theorem \ref{teor:chern-connection}). We also prove that $(M,J,g)$ is a quasi-K\"{a}hler type manifold  if and only if the Chern and the first canonical connection coincide (Proposition \ref{teor:n0=nc}). After that  we study the well adapted connection. First, we introduce a tensor evaluating the first canonical connection over the condition (\ref{eq:welladapted}) which defines the well adapted connection. This tensor, denoted by $\mathcal F(\nabla^0)$, measures in fact the difference between the first canonical and the well adapted connection (Theorem \ref{teor:nabla0=nablaw}). We show that $\mathcal F(\nabla^0)$ can be obtained from the second Nijenhuis tensor if $\alpha\varepsilon=-1$, and from the Nijenhuis tensor in the case $\alpha\varepsilon=1$ (see formulas (\ref{eq:functorial-primeracanonica-1}) and (\ref{eq:functorial-primeracanonica+1})). We continue our study about this connection relating the torsion tensor with the class of the manifold, principally with K\"{a}hler type manifolds. Following the above quoted order, we also study the Kobayashi-Nomizu and Yano connections (see formulas (\ref{eq:leykn}) and \ref{eq:leyyano}).  We prove that they coincide if and only if $J$ is integrable. We also characterize when they are natural connections (Corollaries \ref{teor:KN-natural} and \ref{teor:yano-natural}). It is a remarkable fact about the Kobayashi-Nomizu connection: it is natural if and only if the manifold is of quasi-K\"{a}hler type. Besides, in this case, the Kobayashi-Nomizu and the well adapted connections coincide. This fact supplies an explicit expression of the well adapted connection more handle than  condition  (\ref{eq:welladapted}).  Finally we analyze the existence of natural connections with totally skew-symmetric torsion. We characterize their existence by conditions valid for the four geometries unified under the notion of $(J^2=\pm1)$-metric manifold (Theorem \ref{teor:condicion-sk}). Of course, we show that these conditions are equivalent to previous chacterizations that assure the existence of such connections (Propositions \ref{teor:sk-1-nijenhuis} and \ref{teor:sk+1-quasikahler}). We also obtain simple expressions of connections with totally skew-symmetric torsion, if there exist, by means of the covariant derivative $\nabla^0$, the metric $g$ and the tensor $\nabla^{\mathrm{g}} J$ (see formulas (\ref{eq:sk-1}) and (\ref{eq:sk+1})).

Section 6 is devoted to the study of canonical connections. In the $\alpha \varepsilon =-1$ case, they were introduced in the papers \cite{gauduchon} and \cite{ivanov}. They are generated by the first canonical and the Chern connections, \textit{i.e.}, they form the one-parameter family $\{ (1-t)\nabla ^{0}+t\nabla ^{\mathrm{c}} \colon t \in \mathbb{R}\}$ (see Theorem \ref{teor:familiauniparametrica-1}).  In this theorem we also prove that this family is also generated by the first canonical and the well adapted connections, and thus one can parametrize it as $\{ (1-s)\nabla ^{0}+s\nabla ^{\mathrm{w}}\colon s\in \mathbb{R}\}$. This is very important, because one can define a family of canonical connections in the case $\alpha \varepsilon =1$, generated by the first canonical and the well adapted connections (Proposition \ref{teor:familiauniparametrica+1}). The key of the definition of these families of natural connections is the tensor $\mathcal F(\nabla^0)$ again, which allows to obtain a unified presentation of them (see formulas (\ref{eq:familiauniparametrica-1}) and (\ref{eq:familiauniparametrica+1})). We prove that all distinguished connections studied in the previous section belong to the one-parameter family of canonical connections of the $(J^2=\pm1)$-metric manifold, when they exist and they are natural connections. 

We will consider smooth manifold and operators being of class $C^{\infty }$.

\section{Reducible connections on $(J^{2}=\pm 1)$-manifolds}

A manifold $M$ having a tensor field $J$ of type $(1,1)$ with $J^{2}=\alpha Id$, where $Id$ denotes the identity tensor field and $\alpha \in \{-1,1\}$, is said to be a $(J^{2}=\pm 1)$-manifold or a manifold endowed with an $\alpha$-structure. The corresponding $G$-structure is denoted by $\pi \colon {\mathcal C}_{\alpha} \to M$ and the manifold is said to have a $G_{\alpha }$-structure. We are interested in the study of connections reducible to a $G_{\alpha}$-structure. The following result characterizes these connections, generalizing the well known result for the almost complex manifolds (\cite[Vol. II, Prop. 3.3]{KN}).

\begin{prop}
Let $(M,J)$ be a $(J^{2}=\pm 1)$-manifold.  A linear connection $\Gamma$ on $M$ is reducible to $\pi \colon {\mathcal C}_{\alpha} \to M$ if and only if its covariant derivative  $\nabla$ parallelizes $J$, i.e., $\nabla J=0$, which means $\nabla _{X}JY =J(\nabla _{X}Y),$ for all vector fields $X,Y$ on $M$. In this case, $\nabla$ is said to be natural or adapted to the $\alpha $-structure.
\end{prop}

Following \cite{etayo} and \cite{rod-nurowski}, one can determine the set of  covariant derivatives  adapted to a $G_{\alpha}$-structure.

\begin{lema}
\label{teor:JconexionesS}
Let $(M,J)$ be a manifold endowed with an  $\alpha$-structure $J$ and let $\nabla$ be a  covariant derivative on $M$. The set of covariant derivatives adapted  to $J$ is:
\[
\{ \nabla + S \colon S \in\mathcal T^1_2 (M), (\nabla_X J)Y = JS(X,Y)-S(X,JY), \forall X, Y \in \mathfrak{X} (M) \}.
\]
\end{lema}

In \cite[Lemma 4.3]{racsam} we have obtained the corresponding result to $(J^{2}=\pm 1)$-metric manifolds. The main difference between both results is the existence of the Levi-Civita connection in the metric case, which allows to parametrize all the adapted connections by the potential tensor field (the difference of any connection and the Levi Civita connection). In the case of an $\alpha$-structure we have no such a distinguished connection. This fact implies that one can not define specific connections obtained from a distinguished one. Thus, we can define connection types, but no isolated connections. In this section we deal with the following types of connections on $(M,J)$:

\begin{itemize}
\item $\nabla ^{0}$ type connections are natural connections which can be used to parametrize the set of natural connections taking a natural connection as starting point. Each covariant derivative $\nabla$  on the manifold defines a $\nabla ^{0}$ type connection.

\item $\nabla ^{1}$ type connections  are also  natural connections. The corresponding covariant derivatives will be denoted as $\nabla ^{1}$. Each $\nabla ^{0}$ type connection defines a $\nabla ^{1}$ connection. All the connections $\nabla ^{s}$ of the 1-parameter family of  connections  $\{ \nabla ^{s}=(1-s)\nabla ^{0}+s\nabla ^{1}\colon s\in \mathbb{R}\}$ are natural connections. Each covariant derivative  $\nabla$ defines a unique $\nabla ^{0}$ type connection and then a unique $\nabla ^{1}$ type connection. When $\nabla$ is torsion-free, the corresponding $\nabla ^{1}$ type connection is called a Kobayashi-Nomizu type connection.

\item Yano type connections, whose covariant derivative will be denoted as $\widetilde{\nabla }^1$, are natural if and only if $J$ is integrable. Each torsion-free connection on $M$ induces a Yano type connection.
\end{itemize}

Besides proving the above results we will solve other problems such as  the characterization of the coincidence among $\nabla$, $\nabla ^{0}$, $\nabla^{1}$ and $\widetilde{\nabla}^1$.
\bigskip

In order to have a better presentation of the set of adapted connections, we observe that one can decompose  $\mathcal T^1_2 (M)$ as direct sum of two suitable subspaces:

\begin{prop}
\label{teor:decomposition}
Let $(M,J)$ be a $(J^{2}=\pm 1)$-manifold. Then
$ \mathcal T^1_2 (M)= \mathcal A_{\alpha} \oplus \mathcal L_{\alpha}$,
where
\begin{eqnarray*}
\mathcal A_{\alpha} &=& \{ K   \in \mathcal T^1_2 (M) \colon  K(X,JY)=-JK(X,Y), \forall X, Y \in \mathfrak{X} (M) \},\\
\mathcal L_{\alpha} &=& \{ Q  \in \mathcal T^1_2 (M) \colon Q(X,JY)=JQ(X,Y), \forall X, Y \in \mathfrak{X} (M) \}.
\end{eqnarray*}
\end{prop}

{\bf Proof.}  Let $S \in \mathcal T^1_2 (M)$ and
\[
K(X,Y)= \frac{1}{2} (S(X,Y)-\alpha JS(X,JY)), \quad
Q(X,Y)= \frac{1}{2} (S(X,Y)+\alpha JS(X,JY)),  \quad \forall X, Y \in \mathfrak{X} (M),
\]
then $S=K+Q$. Besides one has
\begin{eqnarray*}
K(X, JY)=\frac{1}{2} (S(X,JY)-JS(X,Y)), &\quad& JK(X,Y)=\frac{1}{2} (JS(X,Y)-S(X,JY)), \\
Q(X, JY)=\frac{1}{2} (S(X,JY)+JS(X,Y)), &\quad& JQ(X,Y)=\frac{1}{2} (JS(X,Y)+S(X,JY)), \quad \forall X, Y \in \mathfrak{X} (M),
\end{eqnarray*}
thus proving $K \in \mathcal A_{\alpha}$ and $Q \in \mathcal L_{\alpha}$. As $\mathcal A_{\alpha}\cap \mathcal L_{\alpha}=\emptyset$, one obtains $ \mathcal T^1_2 (M)= \mathcal A_{\alpha} \oplus \mathcal L_{\alpha}$. $\blacksquare$

\subsection{$\nabla ^{0}$ type connections}

Let $\nabla$ be a covariant derivative on a $(J^2=\pm1)$-manifold  $(M,J)$. Then, we denote by   $\nabla^{0}$ the covariant derivative given by
\begin{equation}
\nabla^{0}_X Y = \nabla_X Y +\frac{(-\alpha)}{2} (\nabla_X J) JY, \quad \forall X, Y \in \mathfrak{X} (M).
\label{eq:Jconexion0}
\end{equation}
One easily checks that $ \nabla^{0}_X (JY) = (1/2) \nabla _{X}JY + (1/2) J(\nabla _{X}Y)= J(\nabla^{0}_X Y)$ thus proving  $\nabla^{0}J=0$. i.e., $\nabla^{0}$ is a natural connection respect to $J$. This connection can be used to parametrize all the adapted covariant derivatives. Observe, nevertheless, that $\nabla^{0}$ is not uniquely defined (depends on the arbitrary derivative $\nabla$).

\begin{prop}
\label{teor:JconexionesQ}
Let $(M,J)$ be a manifold endowed with an $\alpha$-structure  $J$  and let $\nabla$ be a covariant derivative on $M$. The set of derivation laws adapted to $J$ is:
\[
\{ \nabla^{0} + Q \colon Q \in  \mathcal L_{\alpha} \}.
\]
\end{prop}

{\bf Proof.} Let $\nabla^{\mathrm{a}}=\nabla ^{0}+Q$, with $Q\in \mathcal L_{\alpha}$. A direct calculus shows that $\nabla^{\mathrm{a}}J=0$. Conversely let $\nabla^{\mathrm{a}}$ be a natural connection and consider the difference tensor $Q=\nabla^{\mathrm{a}}-\nabla ^{0}$. Then, a direct calculus shows that $Q(X,JY)=JQ(X,Y)$ for all vectors fields $X, Y$ on $M$,  thus proving $Q\in  \mathcal L_{\alpha}$. $\blacksquare$

\begin{obs} The above result shows that one can parametrize the set of natural connections taking the natural connection $\nabla ^{0}$ as a starting point, which differs from the case in Lemma \ref{teor:JconexionesS} where the connection $\nabla$ does not have  to be natural. Besides
\[
K(X,Y)=\frac{(-\alpha)}{2} (\nabla_X J) JY, \quad \forall X, Y \in \mathfrak{X} (M),
\]
is the unique tensor in $\mathcal A_{\alpha} $ such that $\nabla + K$ is a natural connection. This is the key point in the definition of $\nabla ^{0} =\nabla +K$. We prove this claim.  Let $\nabla^{\mathrm{a}}$ be  a natural connection. Taking into account Lemma \ref{teor:JconexionesS} and Proposition \ref{teor:decomposition} we can decompose $\nabla^{\mathrm{a}}=\nabla +S = \nabla + K +Q$ where
\[
 (\nabla_X J)Y = JS(X,Y)-S(X,JY), \quad \forall X, Y \in \mathfrak{X} (M), \quad K \in \mathcal A_{\alpha},  Q \in \mathcal L_{\alpha}.
\]
Then one has
\[(\nabla_X J)Y = JS(X,Y)-S(X,JY) = JK(X,Y)-K(X,JY) + JQ(X,Y)-Q(X,JY)= -2K(X,JY)
\]
and then
\[ K(X,JY) =-\frac{1}{2}(\nabla_X J)Y\quad \Rightarrow \quad \alpha K(X,Y) = K(X,J^{2}Y) = -\frac{1}{2}(\nabla_X J)JY
\]
thus proving $K(X,Y)=\frac{(-\alpha)}{2} (\nabla_X J) JY$.
\end{obs}

\subsection{Kobayashi-Nomizu type connections}

Let $\nabla$ be a covariant derivative on a $(J^2=\pm1)$-manifold $(M,J)$. We have found a natural covariant derivative $\nabla ^{0}$ defined from $\nabla$ as $\nabla ^{0}=\nabla +K$ as it is specified in formula (\ref{eq:Jconexion0}). We are looking for another natural covariant derivative $\nabla ^{1}$ defined from $\nabla$. Let us consider the tensor field  $L_K$, defined as
\[
L_K(X,Y)=\frac{1}{2}(K(Y,X)-\alpha K(JY,JX)) =  \frac{(-\alpha)}{4}((\nabla_YJ) JX - (\nabla_{JY} J) X),  \quad \forall X, Y \in \mathfrak{X} (M).
\]
One easily checks that $L_K \in \mathcal L_{\alpha}$.  Then,
 $\nabla^1 = \nabla ^{0}+L_{K}$ is a natural covariant derivative according to  Proposition \ref{teor:JconexionesQ}. The explicit expression is:
\begin{equation}
\nabla^1_X Y = \nabla^{0}_X Y + \frac{(-\alpha)}{4}((\nabla_YJ) JX - (\nabla_{JY} J) X), \quad   \forall X, Y \in \mathfrak{X} (M).
\label{eq:Jconexion1}
\end{equation}

\begin{prop}
\label{teor:uniparametrica-alpha}
Let $(M,J)$ be a manifold endowed with an $\alpha$-structure  $J$ and let $\nabla$ be a covariant derivative on  $M$. Then
\[
\nabla^{\mathrm{s}}_X Y = (1-s) \nabla^{0}_X Y +  s \nabla^1_X Y, \quad \forall X, Y, \in \mathfrak{X} (M), s \in \mathbb{R},
\]
is a one-parameter family of natural covariant derivatives.  The explicit expression of $\nabla^{\mathrm{s}}$ is given by
\[
\nabla^{\mathrm{s}}_X Y = \nabla^{0}_X Y   + \frac{(-\alpha)s}{4}((\nabla_YJ) JX - (\nabla_{JY} J) X), \quad \forall X, Y \in \mathfrak{X} (M), s \in \mathbb{R}.
\]
\end{prop}

{\bf Proof.} Trivial. $\blacksquare$
\bigskip

Up to this point, we have found a one-parameter family of covariant derivatives adapted to an $\alpha$-structure, which is defined from an arbitrary covariant derivative $\nabla$. The family is determined by $\nabla ^{0}$ and $\nabla ^{1}$.  One expect  better properties of $\nabla ^{0}$ and $\nabla ^{1}$ if $\nabla$ is a torsion-free covariant derivative. This is not a restriction, because, as is well known, one also can define a torsion-free covariant derivative $\tilde{\nabla }$ from $\nabla$, given by $\tilde{\nabla }_{X}Y= \nabla _{X}Y-\frac{1}{2}\mathrm{T}(X,Y)$, for all vectors fields $X, Y$ on $M$, where $\mathrm{T}$ denotes the torsion tensor of $\nabla$. In fact, one has:

\begin{prop}
\label{teor:nosimetrica}
Let $M$ be a manifold and let $\nabla$ be a torsion-free covariant derivative on $M$. The torsion tensor  $\mathrm{T}^S$ of the covariant derivative  $\nabla^S = \nabla+ S$, with $S \in \mathcal T^{1}_2(M)$,  satisfies
\[
\mathrm{T}^S (X,Y)=S(X,Y)-S(Y,X),   \quad \forall X, Y \in \mathfrak{X} (M).
\]
\end{prop}

{\bf Proof.} Given $X, Y$  vector fields on $M$ one has
\[
\mathrm{T}^S (X,Y) = \nabla_X Y + S(X,Y)- \nabla _Y X- S(Y,X)-[X,Y]= S(X,Y)-S(Y,X). \ \blacksquare
\]

Thus, in the rest of this section, we focus our attention on  torsion-free covariant derivatives. Let us remember

\begin{defin}
Let $M$ be a manifold and let $J$ be a tensor field of type  $(1,1)$. The Nijenhuis tensor of $J$ is the tensor field of type  $(1,2)$ given by
\[
N_J(X,Y) = J^2[X,Y]+[JX,JY]-J[JX,Y]-J[X,JY],\quad  \forall X,Y \in \mathfrak{X} (M).
\]
\end{defin}

Then we obtain some easy results which will be used in the future.

\begin{prop}
Let  $(M,J)$ be a $(J^2=\pm 1)$-manifold and let $\nabla$ be a torsion-free covariant derivative on $M$. Then
\begin{equation}
N_J(X,Y) =(\nabla_X J) JY + (\nabla_{JX} J) Y - (\nabla_Y J) JX -(\nabla_{JY} J) X,   \quad \forall X, Y \in \mathfrak{X} (M).
\label{eq:nijenhuis2}
\end{equation}
\end{prop}

{\bf Proof.} Taking into account $\nabla$ is torsion-free one has  $\nabla_X Y -\nabla_Y X = [X,Y]$, $\forall X, Y \in \mathfrak{X} (M)$, the proof is trivial. $\blacksquare$

\begin{prop}
\label{teor:JpropiedadKN}
Let  $(M,J)$ be a $(J^2=\pm 1)$-manifold and let $\nabla$ be a torsion-free covariant derivative on $M$. Then the natural covariant  derivative $\nabla^1$ defined in (\ref{eq:Jconexion1})
satisfies
\[
(-\alpha) N_J(X,Y)=4 \mathrm{T}^{1} (X, Y),   \quad \forall X, Y \in \mathfrak{X} (M),
\]
where $\mathrm{T}^{1}$ denotes the torsion tensor of $\nabla^1$.
\end{prop}

{\bf Proof.} Formula (\ref{eq:Jconexion0}) expresses $\nabla^{0}$ in terms of $\nabla$ and formula (\ref{eq:Jconexion1})  expresses $\nabla^{1}$ in terms of $\nabla^{0}$. Combining both formulas one has
\begin{equation}
\nabla^1_X Y = \nabla_X Y +\frac{(-\alpha)}{2} (\nabla_X J) JY + \frac{(-\alpha)}{4}((\nabla_YJ) JX - (\nabla_{JY} J) X), \quad   \forall X, Y \in \mathfrak{X} (M).
\label{eq:Jconexion1larga}
\end{equation}
Applying Proposition \ref{teor:nosimetrica} to the difference tensor  $S=\nabla ^{1}-\nabla$, one obtains
\begin{eqnarray*}
\mathrm{T}^{1}(X, Y) &=& S(X, Y)-S(Y,X)\\
                                           &=& \frac{(-\alpha)}{2}(\nabla_X J) JY + \frac{(-\alpha)}{4}((\nabla_YJ) JX - (\nabla_{JY} J) X)\\
                                           &-& \frac{(-\alpha)}{2}(\nabla_Y J) JX - \frac{(-\alpha)}{4}((\nabla_X J) JY - (\nabla_{JX} J) Y)\\
                                           &=& \frac{(-\alpha)}{4} ( (\nabla_X J) JY + (\nabla_{JX} J) Y - (\nabla_Y J) JX - (\nabla_{JY} J) X)= \frac{(-\alpha)}{4} N_J(X,Y) ,
\end{eqnarray*}
for all vector fields $X$, $Y$ on $M$, where the last equality follows from (\ref{eq:nijenhuis2}). $\blacksquare$

\begin{obs} In the case $\alpha =-1$ the adapted covariant derivative $\nabla ^{1}$ had been previously studied by Kobayashi and Nomizu (see \cite[Vol. II, Theor. 3.4]{KN}). Starting from a torsion-free covariant derivative $\nabla$, they had introduced the covariant derivative  $\widetilde \nabla$ as
\[
\widetilde \nabla_X Y = \nabla_X Y - Q(X,Y), \quad \forall X, Y \in \mathfrak{X}(M),
\]
where
\[
Q(X,Y)=\frac{1}{4} ((\nabla_{JY} J)X +J((\nabla_Y J)X) + 2J((\nabla_X J) Y)), \quad \forall X, Y \in \mathfrak{X}(M).
\]
Given $X, Y$ vector fields on $M$ and defining
\[
S(X,Y)= -Q(X,Y)= -\frac{1}{2} J((\nabla_X J) Y)  -\frac{1}{4} (J((\nabla_Y J)X) +(\nabla_{JY} J)X ),
\]
by property  (\ref{eq:tensorNJ1}) (which is true for any covariant derivative) the above expression reads as
\[
S(X,Y)= \frac{1}{2} (\nabla_X J) JY  + \frac{1}{4} ( (\nabla_Y J)JX  - (\nabla_{JY} J)X ),
\]
which is formula (\ref{eq:Jconexion1})
\[
S(X, Y)= \frac{(-\alpha)}{2} (\nabla_X J) JY + \frac{(-\alpha)}{4}  ( (\nabla_Y J)JX  - (\nabla_{JY} J)X), \quad \forall X, Y \in \mathfrak{X}(M),
\]
in the case $\alpha=-1$, thus proving $\widetilde \nabla =\nabla^1$.
\end{obs}

The above Remark allows us to introduce the following:

\begin{defin}
Let  $(M,J)$ be a $(J^2=\pm 1)$-manifold and let $\nabla$ be a torsion-free covariant derivative on $M$. The adapted covariant derivative $\nabla^1$ on $M$ defined from $\nabla$ in (\ref{eq:Jconexion1}) is said to be a covariant derivative of Kobayashi-Nomizu type.
\end{defin}

One can prove the following result which shows the interest of derivatives of Kobayashi-Nomizu type:

\begin{prop}
\label{teor:leyKN-integrable}
Let $(M,J)$ be a manifold endowed with an $\alpha$-structure $J$. The following conditions are equivalent:
\begin{enumerate}
\renewcommand*{\theenumi}{\roman{enumi})}
\renewcommand*{\labelenumi}{\theenumi}

\item The $\alpha$-structure $J$ is integrable.

\item The manifold $M$ admits a torsion-free covariant derivative adapted to $J$.
\end{enumerate}
\end{prop}

{\bf Proof.}

$i) \Rightarrow ii)$ As $J$ is integrable, then $N_{J}=0$. Taking into account Proposition \ref{teor:JpropiedadKN} one obtains that any covariant of  Kobayashi-Nomizu type is torsion-free and natural.

$ii) \Rightarrow i)$ Let $\nabla$ be a torsion-free adapted covariant derivative. Then, by formula  (\ref{eq:nijenhuis2}),  one obtains that the Nijenhuis tensor $N_J$ vanishes and thus $J$ is integrable. $\blacksquare$

\begin{obs} Let $(M,J)$ be a manifold endowed with an $\alpha$-structure $J$. If $\nabla$ is a torsion-free adapted covariant derivative, then $\nabla =\nabla ^{0}=\nabla ^{1}$ because of formulas (\ref{eq:Jconexion0})  and (\ref{eq:Jconexion1larga}).
\end{obs}

\subsection{Yano type connections}

A Yano type connection is defined from a torsion-free connection as follows:

\begin{defin}
Let  $(M,J)$ be a $(J^2=\pm 1)$-manifold and let $\nabla$ be a torsion-free covariant derivative on $M$. The covariant derivative $\widetilde\nabla^1$ on $M$ defined from $\nabla$ as  follows
\begin{equation}
\widetilde\nabla^1_X Y = \nabla_X Y + \frac{(-\alpha)}{2} (\nabla_Y J) JX +\frac{(-\alpha)}{4} ((\nabla_X J)JY - (\nabla_{JX} J)Y), \quad \forall X, Y \in \mathfrak{X} (M),
\label{eq:leyyano2}
\end{equation}
is said to be a covariant derivative of Yano type.
\end{defin}

Yano had defined special connections in both the almost complex and almost product cases, which are the model for the above definition. For instance, in \cite{yano} he studied an almost product manifold $(M,J)$ and by means of the Levi Civita connection of an arbitrary metric $g$, he defined
\[
\nabla_X Y = \nabla^{\mathrm{g}}_X Y - \frac{1}{2} (\nabla^{\mathrm{g}}_Y J) X -\frac{1}{4} ((\nabla^{\mathrm{g}}_X J)JY - (\nabla^{\mathrm{g}}_{JX} J)Y), \quad \forall X, Y \in \mathfrak{X} (M).
\]

This was important because of the following result:

\begin{teor}[{\cite[Theor. 25]{yano}}]
Let $(M,J)$ be an almost product manifold. The following conditions are equivalent:
\begin{enumerate}
\renewcommand*{\theenumi}{\roman{enumi})}
\renewcommand*{\labelenumi}{\theenumi}

\item The almost product structure $J$ is integrable.

\item The manifold $M$ admits a torsion-free connection adapted to $J$.
\end{enumerate}
\end{teor}

In order to prove the result, one needs $\nabla^{\mathrm{g}}$ to be torsion-free. In fact, this is the  essential point, and not other properties of $\nabla^{\mathrm{g}}$. Thus, the above result remains true when $\nabla^{\mathrm{g}}$ is substituted by any torsion-free connection. And this is the reason of our above definition. That definition is quite similar to that of Kobayashi-Nomizu type connections given in (\ref{eq:Jconexion1larga}), thus leading us to study the relationship between these Yano and Kobayashi-Nomizu types connections in the case they are derived from the same torsion-free connection. We need the following technical lemma in order to answer the question.

\begin{lema}
\label{teor:yanotecnico}
Let  $(M,J)$ be a $(J^2=\pm 1)$-manifold and let $\nabla$ be a torsion-free covariant derivative on $M$. Let $\nabla^1$ and $\widetilde \nabla^1$ be the  Kobayashi-Nomizu and the Yano type connections defined from  $\nabla$ according to  (\ref{eq:Jconexion1larga}) and (\ref{eq:leyyano2}) respectively.  For all vector fields $X, Y$ on $M$ the following relations hold:
\begin{enumerate}
\renewcommand*{\theenumi}{\roman{enumi})}
\renewcommand*{\labelenumi}{\theenumi}

\item $\nabla^1_X Y - \widetilde\nabla^1_X Y  = \frac{(-\alpha)}{4} N_J(X,Y)$.

\item $\widetilde{\mathrm{T}}^1 (X,Y)= \frac{\alpha}{4} N_J(X,Y)$, where $\widetilde{\mathrm{T}}^1$ denotes the torsion tensor of $\widetilde\nabla^1$.

\item $\nabla^1_X Y= \widetilde\nabla^1_X Y- \widetilde{\mathrm{T}}^1 (X,Y)$.

\item $J\widetilde{S}^1(X,Y)-\widetilde{S}^1(X,JY)-(\nabla_X J) Y = \frac{(-\alpha)}{2} N_J(JX,Y)$, where $\widetilde{S}^1=\widetilde\nabla^1-\nabla$.
\end{enumerate}
\end{lema}

{\bf Proof.}

$i)$ According to formulas (\ref{eq:Jconexion1larga}), (\ref{eq:leyyano2}) and (\ref{eq:nijenhuis2}) one has
\begin{eqnarray*}
\nabla^1_X Y  - \widetilde\nabla^1_X Y  &=& \frac{(-\alpha)}{2} (\nabla_X J) JY + \frac{(-\alpha)}{4}((\nabla_YJ) JX - (\nabla_{JY} J) X) \\
                               &+& \frac{\alpha}{2} (\nabla_Y J) JX +  \frac{\alpha}{4}((\nabla_X J) JY - (\nabla_{JX} J) Y)\\
                               &=& \frac{(-\alpha)}{4} ((\nabla_X J) JY + (\nabla_{JX} J) Y - (\nabla_Y J) JX -(\nabla_{JY} J) X) = \frac{(-\alpha)}{4} (N_J(X,Y), \quad \forall X, Y \in \mathfrak{X} (M).
\end{eqnarray*}

$ii)$ As $\nabla$ is torsion-free, according to Proposition \ref{teor:nosimetrica} and formula (\ref{eq:nijenhuis2}) one obtains
\begin{eqnarray*}
\widetilde{\mathrm{T}}^1 (X, Y) &=&  \widetilde{S}^1(X,Y)- \widetilde{S}^1(Y,X)\\
                               &=& \frac{(-\alpha)}{2} (\nabla_Y J) JX +  \frac{(-\alpha)}{4}((\nabla_X J) JY - (\nabla_{JX} J) Y)\\
                               &+&\frac{\alpha}{2} (\nabla_X J) JY + \frac{\alpha}{4}((\nabla_YJ) JX - (\nabla_{JY} J) X) \\
                               &=& \frac{\alpha}{4} ((\nabla_X J) JY + (\nabla_{JX} J) Y - (\nabla_Y J) JX -(\nabla_{JY} J) X) = \frac{\alpha}{4} (N_J(X,Y), \quad \forall X, Y \in \mathfrak{X} (M).
\end{eqnarray*}

$iii)$ Trivial, by the previous items $i)$ and $ii)$.

$iv)$ Given $X, Y$ vector fields on $M$ and taking into account formulas (\ref{eq:leyyano2}) and (\ref{eq:tensorNJ1}) one has
\begin{eqnarray*}
J\widetilde{S}^1(X,Y)&=& \frac{1}{2} (\nabla_Y J) X +  \frac{1}{4}(\nabla_X J) Y -\frac{\alpha}{4} (\nabla_{JX} J) JY,\\
\widetilde{S}^1(X,JY)&=&\frac{(-\alpha)}{2} (\nabla_{JY} J) JX -  \frac{1}{4}((\nabla_X J) Y+ \frac{\alpha}{4} (\nabla_{JX} J) JY,
\end{eqnarray*}
and, according to formula (\ref{eq:nijenhuis2}) one obtains
\[
J\widetilde{S}^1(X,Y)-\widetilde{S}^1(X,JY) -(\nabla_X J) Y = -\frac{1}{2}(\nabla_X J) Y -\frac{\alpha}{2} (\nabla_{JX} J) JY + \frac{1}{2} (\nabla_Y J) X +\frac{\alpha}{2} (\nabla_{JY} J) JX = \frac{(-\alpha)}{2} N_J(JX, Y). \  \blacksquare
\]

Then, the above Lemma and  Lemma \ref{teor:JconexionesS} allow to obtain:

\begin{prop}
\label{Yanoadaptada}
Let  $M$ be a $(J^2=\pm 1)$-manifold and let $\nabla$ be a torsion-free covariant derivative on $M$. Let $\nabla^1$ and $\widetilde \nabla^1$ be the  Kobayashi-Nomizu and the Yano type connections defined from  $\nabla$ according to  (\ref{eq:Jconexion1larga}) and (\ref{eq:leyyano2}) respectively.  Then the following conditions are equivalent:
\begin{enumerate}
\renewcommand*{\theenumi}{\roman{enumi})}
\renewcommand*{\labelenumi}{\theenumi}

\item The $\alpha$-structure $J$ is integrable.

\item The covariant derivative $\widetilde\nabla^1$ is adapted to  $J$.

\item The covariant derivatives $\widetilde\nabla^1$ and $\nabla^1$ coincide.
\end{enumerate}
\end{prop}

As one can see, Kobayashi-Nomizu and Yano type connections derived from the same torsion-free connection are almost equal and then the above result can be derived from Proposition \ref{teor:leyKN-integrable}. In fact, naming $\nabla ^{1}=\nabla +S$ and $\widetilde{\nabla }^1= \nabla +\widetilde{S}^1$ one has
\[
S^1(X,Y)=\widetilde{S}^1(Y,X), \quad \forall X, Y \in \mathfrak{X} (M),
\]
and then their torsion tensors satisfy $\mathrm{T^1}(X,Y) = -\widetilde{\mathrm{T}}^1(X,Y)$. Then, why is the interest in having these two connection types? We will study the question on $(J^{2}\pm 1)$-metric manifolds, where we have the Levi Civita connection $\nabla^{\mathrm{g}}$ as a distinguished torsion-free covariant derivative. Then Kobayashi-Nomizu and Yano connections are uniquely determined. The first one  and the well adapted connection coincide in the case of quasi-K\"{a}hler manifolds while the Yano connection is torsion-free if and only if the structure $J$ is integrable. All of this will be showed later.

\section{Reducible connections on  $(J^2=\pm1)$-metric manifolds}

The core of this paper concerns to manifolds having two compatible structures, an $\alpha$-structure $J$ and a (semi)-Riemannian metric $g$. Compatibility means:

\begin{defin}[{\cite[Defin. 3.1]{racsam}}]
Let $M$ be a  manifold, $g$ a semi-Riemannian metric on $M$, $J$ a  tensor field of type (1,1) and $\alpha ,\varepsilon \in \{-1,1\}$. Then $(J,g)$ is called an $(\alpha ,\varepsilon )$-structure on $M$ if
\[
J^{2} = \alpha Id, \quad \mathrm{trace} \, J=0, \quad g(JX,JY)= \varepsilon g(X,Y), \quad \forall X, Y \in \mathfrak{X}(M), \]
$g$ being a Riemanianan metric if $\varepsilon =1$. Then  $(M,J,g)$ is called a $(J^{2}=\pm 1)$-metric manifold.
\end{defin}

Condition $\mathrm{trace} \, J=0$ is a consequence of the other conditions in all the cases unless the $(1,1)$. We impose it in this case looking for a common treatment of all the four geometric structures. See  \cite{racsam} for a more complete description. Having a metric we can choose its Levi Civita connection as the starting point to study connections on a $(J^{2}=\pm 1)$-metric manifold. This is the key point in which manifolds endowed with an $\alpha$- or an $(\alpha,\varepsilon )$-structure differ.
\bigskip

The $G$-structure defined by an  $(\alpha ,\varepsilon )$-structure will be denoted as a $G_{ (\alpha ,\varepsilon )}$-structure. The corresponding structure groups and Lie algebras have been studied in \cite{racsam}. In particular one has:

\begin{prop}[{\cite[Prop. 3.9]{racsam}}]
Let $(M,J,g)$ be a  $(J^2=\pm1)$-metric manifold and let $\pi \colon {\mathcal C}_{(\alpha ,\varepsilon )} \to M$ the  $G_{(\alpha ,\varepsilon )}$-structure on $M$  defined by $(J,g)$. Let $\Gamma$ be a linear connection on  $M$ and let $\nabla$ be the corresponding derivation law. Then $\Gamma$ is a reducible connection to $\pi \colon {\mathcal C}_{(\alpha ,\varepsilon )} \to M$ if and only if $\nabla J = 0, \nabla g = 0$.
\end{prop}

As in the case of an $\alpha$-structure, we introduce the following:

\begin{defin}[{\cite[Def. 4.1]{racsam}}]
Let $(M,J,g)$ be a $(J^2=\pm1)$-metric manifold.  A  covariant derivative or derivation law $\nabla^{\mathrm{a}}$ on $M$ is said to be natural or adapted to $(J,g)$ if $\nabla^{\mathrm{a}} J=0, \nabla^{\mathrm{a}} g=0$.
\end{defin}

As we have a distinguished derivative, that defined by the Levi Civita connection, we can compare any other one with that one:

\begin{defin}[{\cite[Prop. 4.2]{racsam}}]
\label{teor:tensorpotencial}
Let $(M,J,g)$ be a $(J^2=\pm1)$-metric manifold, let  $\nabla^{\mathrm{g}}$ be the derivation law of the Levi Civita connection of $g$ and let $\nabla^{\mathrm{a}}$ be a derivation law adapted to $(J,g)$. The potential tensor of $\nabla^{\mathrm{a}}$ is the  tensor $S\in \mathcal T^1_2 (M)$ defined as
\[
S(X,Y)=\nabla^{\mathrm{a}}_X Y -\nabla^{\mathrm{g}}_X Y,   \quad \forall X, Y \in \mathfrak{X} (M).
\]
\end{defin}

Then, we can parametrize the set of natural covariant derivatives  by means of the Levi Civita connection and the potential tensor:

\begin{lema}[{\cite[Lemma 4.3]{racsam}}]
\label{teor:natural}
Let $(M,J,g)$ be a $(J^2=\pm1)$-metric manifold. The set of derivation laws adapted to $(J,g)$ is:
\[
\left\{ \nabla^{\mathrm{g}}+S \colon S \in \mathcal T^1_2 (M),
\begin{array}{l}
JS(X,Y)-S(X,JY)=(\nabla^{\mathrm{g}}_X J) Y, \\
 g(S(X,Y),Z)+g(S(X,Z),Y)=0,
 \end{array}
 \forall X, Y, Z \in \mathfrak{X} (M) \right\}.
\]
\end{lema}

The following result can be thought as a translation of Proposition \ref{teor:nosimetrica} to the present situation:

\begin{prop} Let $(M,J,g)$ be a $(J^2=\pm1)$-metric manifold, let $\nabla^{\mathrm{a}}$ be an adapted covariant derivative and let $\mathrm{T}^{\mathrm{a}}$  (resp. $S$) be the torsion tensor (resp. the potential tensor) of $\nabla^{\mathrm{a}}$. The following equalities hold:
\begin{eqnarray}
\mathrm{T}^{\mathrm{a}} (X,Y)&=&S(X,Y)-S(Y,X),   \quad \forall X, Y \in \mathfrak{X} (M).
\label{teor:tora-potencial}\\
g(S(X,Y),Z)&=& \frac{1}{2} (g(\mathrm{T}^{\mathrm{a}}(X,Y),Z)-g(\mathrm{T}^{\mathrm{a}}(Y,Z),X)+g(\mathrm{T}^{\mathrm{a}}(Z,X),Y)),  \quad  \forall X, Y, Z \in \mathfrak{X} (M).
\label{eq:potential-torsion}
\end{eqnarray}
\end{prop}

{\bf Proof.} Formula (\ref{teor:tora-potencial}) is a direct consequence of Proposition \ref{teor:nosimetrica}. We prove the other formula. As  $\nabla^{\mathrm{a}}$ is adapted to $(J,g)$, according to Lemma \ref{teor:natural}, one has
\begin{eqnarray*}
(g(S(X,Y),Z)+g(S(X,Z,Y))) - (g(S(Z,X),Y)+g(S(Z,Y),X)) &+& \\
(g(S(Y,Z),X)+g(S(Y,X),Z)) + (g(S(X,Y),Z)-g(S(X,Y),Z)) &=& 0, \quad \forall X, Y, Z \in \mathfrak{X} (M),
\end{eqnarray*}
which, taking into account formula  (\ref{teor:tora-potencial}),  reads as
\[
2 g(S(X,Y),Z)+ g(\mathrm{T}^{\mathrm{a}}(X,Z),Y)+g(\mathrm{T}^{\mathrm{a}}(Y,Z),X)+g(\mathrm{T}^{\mathrm{a}}(Y,X),Z)=0. \ \blacksquare
\]

The above properties are well known. Formula (\ref{eq:potential-torsion}) in the Riemannian case appears in
 \cite[Prop. 2.1]{agricola} and \cite[Theor. 3.4]{rod-nurowski}. Both formulas are also used in  \cite{ganchev-mihova} and \cite{mihova}.
\bigskip

The following two results summarize some properties which have easy proofs.

\begin{lema}
\label{teor:primerapropiedad}
 Let $M$ be a manifold endowed with an $\alpha$-structure $J$ and a metric $g$. The following  conditions are equivalent:
\begin{enumerate}
\renewcommand*{\theenumi}{\roman{enumi})}
\renewcommand*{\labelenumi}{\theenumi}

\item $g(JX,JY)=\varepsilon g(X,Y)$, for all  vector fields  $X, Y$ on $M$.

\item $g(JX, Y)=\alpha \varepsilon g(X, JY)$, for all  vector fields  $X, Y$ on $M$.
\end{enumerate}
\end{lema}

\begin{lema}
\label{teor:tensorNJ}
Let $(M,J,g)$ be a  $(J^2=\pm 1)$-metric manifold. The tensor $\nabla^{\mathrm{g}} J$ satisfies the following relations:
\begin{eqnarray}
(\nabla^{\mathrm{g}}_X J) JY &=& - J (\nabla^{\mathrm{g}}_X J) Y, \label{eq:tensorNJ1}\\
g((\nabla^{\mathrm{g}}_X J)Y,Z)&=&\alpha \varepsilon g ((\nabla^{\mathrm{g}}_X J)Z,Y), \label{eq:tensorNJ2}\\
g((\nabla^{\mathrm{g}}_X J) JY, Z) &=& -\alpha \varepsilon g ((\nabla^{\mathrm{g}}_X J) Y, JZ), \label{eq:tensorNJ3}\\
g((\nabla^{\mathrm{g}}_X J) JY, Z) &=& - g((\nabla^{\mathrm{g}}_X J) JZ, Y), \label{eq:tensorNJ4}
\end{eqnarray}
for all  vector fields  $X,Y,Z$ on $M$.
\end{lema}

Following the ideas of the above section about adapted covariant derivatives to an $\alpha$-structure, we can  define $\nabla ^{0}$ as in (\ref{eq:Jconexion0}), choosing $\nabla^{\mathrm{g}}$ as a starting point. Observe that in the present case $\nabla ^{0}$ is uniquely defined on the manifold $M$ because the Levi Civita connection $\nabla^{\mathrm{g}}$  is uniquely determined.

\begin{defin}
\label{teor:first-connection}
Let $(M,J,g)$ be a $(J^2=\pm 1)$-metric manifold. The first canonical connection of $(M,J,g)$ is the linear connection having the covariant derivative $\nabla ^{0}$ given by
\[
\nabla^{0}_X Y= \nabla^{\mathrm{g}}_X Y +\frac{(-\alpha)}{2} (\nabla^{\mathrm{g}}_X J)JY,   \quad \forall X,Y \in \mathfrak{X} (M).
\]
\end{defin}

As one can expect, we have the following result:

\begin{lema}
\label{teor:nabla0natural}
 Let $(M,J,g)$ be a $(J^2=\pm 1)$-metric manifold. Then the covariant derivative $\nabla^{0}$ is adapted to $(J,g)$.
\end{lema}

{\bf Proof.}   According to properties  (\ref{eq:tensorNJ1}) and (\ref{eq:tensorNJ4}), the potential tensor  $S$ of $\nabla^{0}$ satisfies
\begin{eqnarray*}
JS(X,Y)-S(X,JY)&=& \frac{(-\alpha)}{2}J (\nabla^{\mathrm{g}}_X J)JY  +\frac{1}{2}(\nabla^{\mathrm{g}}_X J)Y  = (\nabla^{\mathrm{g}}_X J) Y, \\
g(S(X,Y),Z)+g(S(X,Z),Y)&=& \frac{(-\alpha)}{2} (g((\nabla^{\mathrm{g}}_X J) JY, Z)+g((\nabla^{\mathrm{g}}_X J)JZ,Y)=0,   \quad \forall X, Y, Z \in \mathfrak{X} (M),
\end{eqnarray*}
then by Lemma \ref{teor:natural},  $\nabla^{0}$ is adapted to $(J,g)$. $\blacksquare$
\bigskip

The first canonical connection can be characterized as it is shown in the next proposition, which generalizes that of
 \cite[Theor. 3.4]{rod-nurowski} obtained for the almost Hermitian case.

\begin{prop} Let $(M,J,g)$  be a $(J^2=\pm 1)$-metric manifold. The covariant derivative of the first canonical connection of $(M,J,g)$ is the unique adapted covariant derivative whose potential tensor $S$ satisfies
\[
S(X,Y)+\alpha J S(X,JY)=0, \quad \forall X, Y \in \mathfrak{X}(M).
\]
\end{prop}

{\bf Proof.}  Let $\nabla^{\mathrm{a}}=\nabla^{\mathrm{g}}+S$ be an adapted covariant derivative. Then, for all vector fields $X, Y$ on $M$, one has
\begin{eqnarray*}
S(X,Y)+\alpha J S(X,JY)= 0&\Leftrightarrow& J(S(X,Y)+\alpha J S(X,JY))=0     \\
 &\Leftrightarrow &JS(X,Y)+ S(X,JY)=0.
\end{eqnarray*}
As $\nabla^{\mathrm{a}}$ is an adapted covariant derivative, we have by Lemma \ref{teor:natural}
\[
 JS(X,Y)- S(X,JY)=(\nabla^{\mathrm{g}}_X J) Y,  \quad \forall X, Y \in \mathfrak{X}(M),
\]
and substituting in the above expression one has
\[
S(X,Y)+\alpha J S(X,JY)= 0\Leftrightarrow
2 JS(X,Y)= (\nabla^{\mathrm{g}}_X J) Y \Leftrightarrow S(X,Y)= \frac{(-\alpha)}{2} (\nabla^{\mathrm{g}}_X J) Y \Leftrightarrow \nabla^{\mathrm{a}}=\nabla^{0},
\]
thus proving the result. $\blacksquare$
\bigskip

In Lemma \ref{teor:natural} we have determined the set of adapted covariant derivatives taking the Levi Civita connection of $g$ as the starting point. We can also obtain a result similar to  Proposition \ref{teor:JconexionesQ}, which allows to parametrize that set taking the first canonical connection as starting point.

\begin{lema}
\label{teor:naturalQ}
Let $(M,J,g)$ be a $(J^2=\pm 1)$-metric manifold. The set of natural derivation laws of $(J,g)$ is:
\[
\left\{ \nabla^{0}+Q \colon Q \in \mathcal L_{\alpha}\subseteq \mathcal T^1_2 (M),
 g(Q(X,Y),Z)+g(Q(X,Z),Y)=0, \forall X, Y, Z \in \mathfrak{X} (M) \right\}.
\]
The tensor $Q$ of the natural covariant derivative is said to be the canonical potential tensor.
\end{lema}

{\bf Proof.} Let $\nabla^{\mathrm{a}}=\nabla^{0}+Q$ be a natural covariant derivative, with $Q\in \mathcal T^1_2 (M)$. According to Proposition \ref{teor:JconexionesQ} one knows $Q \in \mathcal L_{\alpha}$ if and only if $\nabla^{\mathrm{a}} J=0$. Thus, we must prove $\nabla^{\mathrm{a}}g=0$.

As $\nabla^{0} g=0$, one has for all vector fields $X, Y, Z$ on $M$
\begin{eqnarray*}
(\nabla^{\mathrm{a}}_X g) (Y,Z) &=& -g(\nabla^{\mathrm{a}}_X Y, Z)-g(\nabla^{\mathrm{a}}_X Z, Y)+ g(\nabla^{0}_X Y, Z)+g(\nabla^{0}_X Z, Y)\\
                         &=&-(g(Q(X,Y),Z)+g(Q(X,Z),Y),
\end{eqnarray*}
thus proving $\nabla^{\mathrm{a}}g=0$ if and only if $g(Q(X,Y),Z)+g(Q(X,Z),Y)=0$. $\blacksquare$
 \bigskip

We have studied the first canonical connection of an $(\alpha , \varepsilon)$-structure  taking in mind the case of  the derivatives $\nabla ^{0}$ associated to an $\alpha$-structure. What can we say about Kobayashi-Nomizu and Yano type covariant derivatives defined in (\ref{eq:Jconexion1larga}) and (\ref{eq:leyyano2})? In Section 5 we will show that in general  they are not reducible connections on $(J^2=\pm1)$-metric manifolds.

\section{Tensors on $(J^2=\pm1)$-metric manifolds defined from $\nabla^{\mathrm{g}} J$}

This is a technical section. We study some tensors derived from $\nabla^{\mathrm{g}} J$ which will be useful in the study of connections in the remaining sections. Given a $(J^2=\pm1)$-metric manifold, we will study three tensors: the covariant derivative $\nabla^{\mathrm{g}} \Phi$ of  the fundamental tensor $\Phi$, the Nijenhuis tensor $N_{J}$ and the second Nijenhuis tensor $\widetilde N_J^{\alpha\varepsilon}$. Properties of the two first tensors are well known and will be summarized. We focus on the expression of all of these tensors by means of $\nabla^{\mathrm{g}} J$.

\subsection{The covariant derivative $\nabla^{\mathrm{g}} \Phi$ of  the fundamental tensor $\Phi$}

Remember the definition:
\begin{defin}
Let $(M,J,g)$ be a $(J^2=\pm 1)$-metric manifold. The fundamental tensor $\Phi$ is the tensor field of type $(0,2)$ defined as
\[
\Phi(X,Y)= g(JX,Y),   \quad \forall X, Y \in \mathfrak{X} (M).
\]
\end{defin}

As is well known, one has:

\begin{lema}
Let $(M,J,g)$ be a $(J^2=\pm 1)$-metric manifold with fundamental tensor $\Phi$.
\begin{enumerate}
\renewcommand*{\theenumi}{\roman{enumi})}
\renewcommand*{\labelenumi}{\theenumi}

\item If $\alpha \varepsilon =-1$ then $\Phi$ is a $2$-form on $M$. In this case $\Phi$ is called the fundamental form of $(M,J,g)$.

\item If $\alpha \varepsilon =1$ then $\Phi$ is a symmetric tensor field. In this case $\Phi=\widetilde g$ is called the twin metric of $g$.
\end{enumerate}
\end{lema}

One can obtain an expression of the covariant derivative of the fundamental tensor by means of $\nabla^{\mathrm{g}} J$:

\begin{prop}
Let $(M,J,g)$ be a $(J^2=\pm 1)$-metric manifold. Then
\begin{equation}
(\nabla^{\mathrm{g}}_X \Phi) (Y,Z)= g((\nabla^{\mathrm{g}}_X J)Y,Z),   \quad \forall X,Y, Z \in \mathfrak{X} (M).
\label{eq:pp-tensorfundamental}
\end{equation}
\end{prop}

{\bf Proof.} Trivial. $\blacksquare$
\bigskip

The above result allows to introduce the more distinguished class of $(J^2=\pm1)$-metric manifolds:

\begin{defin}
Let $(M,J,g)$ be a  $(J^2=\pm1)$-metric manifold. It is said to be a manifold of  K\"{a}hler type if
$\nabla^{\mathrm{g}} \Phi=0.$
\end{defin}

As is well known K\"{a}hler type manifolds are characterized by the condition $\nabla^{\mathrm{g}} J=0$. This condition can be expressed in the following terms:

\begin{lema}
\label{teor:kahler-caracterizacion}
Let $(M,J,g)$ be a $(J^2=\pm1)$-metric manifold. Then $(M,J,g)$ is a K\"{a}hler type manifold if and only if $\nabla^{\mathrm{g}}$ is a covariant derivative adapted to $(J,g)$.
\end{lema}

{\bf Proof.} Trivial, according to formula (\ref{eq:pp-tensorfundamental}). $\blacksquare$
\bigskip

We will end this study of K\"{a}hler type with the following technical result in the case of $(J^2=\pm1)$-metric manifolds with $\alpha \varepsilon =1$.

\begin{lema}
\label{teor:kahler+1tecnico}
Let $(M,J,g)$ be a $(J^2=\pm 1)$-metric manifold with $\alpha \varepsilon =1$. The following conditions are equivalent:
\begin{enumerate}
\renewcommand*{\theenumi}{\roman{enumi})}
\renewcommand*{\labelenumi}{\theenumi}

\item $(M,J,g)$ is a K\"{a}hler type manifold.

\item $g((\nabla^{\mathrm{g}}_X J)Y, Z)= g((\nabla^{\mathrm{g}}_Y J)Z, X) + g((\nabla^{\mathrm{g}}_Z J)X, Y)$, for all vector fields $X, Y, Z$ on $M$.

\item $g((\nabla^{\mathrm{g}}_X J)Y, Y)=0$, for all vector fields $X, Y$ on $M$.
\end{enumerate}
\end{lema}

{\bf Proof.}

$i)\Rightarrow ii)$ Trivial by the definition of a K\"{a}hler type manifold.

$ii)\Rightarrow iii)$ Given $X, Y, Z$  vector fields on $M$ such that $X=Z$  one has
\[
g((\nabla^{\mathrm{g}}_X J)Y, X)= g((\nabla^{\mathrm{g}}_Y J)X, X) + g((\nabla^{\mathrm{g}}_X J)X, Y),
\]
according to formula (\ref{eq:tensorNJ2}) in the case $\alpha \varepsilon =1$ one can deduce $g((\nabla^{\mathrm{g}}_Y J)X, X)=0$.

$iii)\Rightarrow i)$ Given $X, Y, Z$  vector fields on $M$, by property (\ref{eq:tensorNJ2}) in the case $\alpha \varepsilon =1$ one has
\[
g((\nabla^{\mathrm{g}}_X J)Y+Z, Y+Z)= g((\nabla^{\mathrm{g}}_X J)Y, Y)+g((\nabla^{\mathrm{g}}_X J)Y, Z)+g((\nabla^{\mathrm{g}}_X J)Z, Y)+g((\nabla^{\mathrm{g}}_X J)Z, Z)= 2 g(\nabla^{\mathrm{g}}_X J)(Y,Z),
\]
then $\nabla^{\mathrm{g}} \Phi=0$ and thus proving $(M,J,g)$ is a  K\"{a}hler type manifold. $\blacksquare$

\subsection{The Nijenhuis tensor}

According to formula (\ref{eq:nijenhuis2}) the Nijenhuis tensor $N_{J}$ on a  $(J^2=\pm1)$-metric manifold can be written by means of $\nabla^{\mathrm{g}}$ as
\begin{equation}
N_J(X,Y) =(\nabla^{\mathrm{g}}_X J) JY + (\nabla^{\mathrm{g}}_{JX} J) Y - (\nabla^{\mathrm{g}}_Y J) JX -(\nabla^{\mathrm{g}}_{JY} J) X,   \quad \forall X, Y \in \mathfrak{X} (M).
\label{eq:nijenhuis}
\end{equation}
Then one can easily deduce the following result:

\begin{lema}
\label{teor:propiedadesNJ}
Let $(M,J,g)$ be a $(J^2=\pm 1)$-metric manifold. The following relations hold:

\begin{enumerate}
\renewcommand*{\theenumi}{\roman{enumi})}
\renewcommand*{\labelenumi}{\theenumi}
\item
$N_J(Y,X)=-N_J(X,Y), N_J(JX,JY)=\alpha N_J(X,Y), N_J(JX,Y)=N_J(X,JY)$, for all vector fields $X, Y$ on $M$.

\item $g(N_J(JX,Y),JZ)=-\varepsilon g(N_J(X,Y),Z)$, for all vector fields $X, Y, Z$ on $M$.
\end{enumerate}
\end{lema}

The vanishing of the Nijenhuis tensor means the integrability of the $\alpha$-structure $J$. In the case of  $(J^2=\pm1)$-metric manifolds integrability can be expressed in different ways, as the following results show. There exists a difference between cases $\alpha \varepsilon =-1$ and $\alpha \varepsilon =1$ as we are going to show.

\begin{lema}[{\cite[Table I]{gray-hervella}, \cite[Prop.\ 2.1]{olszak}}]
\label{teor:caracterizacionintegrabilidad-1}
Let $(M,J,g)$ be a $(J^2=\pm 1)$-metric manifold with $\alpha \varepsilon =-1$. The following conditions are equivalent:
\begin{enumerate}
\renewcommand*{\theenumi}{\roman{enumi})}
\renewcommand*{\labelenumi}{\theenumi}
\item The Nijenhuis tensor of $J$ vanishes.

\item $(\nabla^{\mathrm{g}}_X J) Y +\alpha (\nabla^{\mathrm{g}}_{JX} J) JY=0$, for all vector fields $X,Y$ on $M$.

\item $(\nabla^{\mathrm{g}}_X \Phi) (Y,Z)+\alpha (\nabla^{\mathrm{g}}_{JX} \Phi) (JY,Z)=0$, for all vector fields $X,Y, Z$ on $M$.\end{enumerate}
\end{lema}

\begin{lema}[{\cite{staikova}, \cite{teofilova}}]
\label{teor:caracterizacionintegrabilidad+1}
Let $(M,J,g)$ be a $(J^2=\pm 1)$-metric manifold with $\alpha \varepsilon =1$. The following conditions are equivalent:
\begin{enumerate}
\renewcommand*{\theenumi}{\roman{enumi})}
\renewcommand*{\labelenumi}{\theenumi}
\item The Nijenhuis tensor vanishes.

\item $g((\nabla^{\mathrm{g}}_{X} J)Y,JZ)+g((\nabla^{\mathrm{g}}_{Y} J)Z,JX)+g((\nabla^{\mathrm{g}}_{Z} J) X, JY)=0$, for all vector fields $X,Y,Z$ on $M$.

\item $(\nabla^{\mathrm{g}}_{X} \Phi) (Y,JZ)+(\nabla^{\mathrm{g}}_{Y} \Phi) (Z,JX)+(\nabla^{\mathrm{g}}_{Z} \Phi) (X, JY)=0$, for all vector fields $X,Y,Z$ on $M$.
\end{enumerate}
\end{lema}

In fact, as we have pointed out, the results have been independently proved for each one of the four geometries. Unifying the study of the four geometries as possible is one of the goals of the present paper.
\bigskip

In order to obtain results about the integrability of $J$ one can also consider the torsion tensor of any  covariant derivative adapted to the $(\alpha,\varepsilon)$-structure. The following two lemmas are examples of this situation.

\begin{lema}
\label{teor:nijenhuis-torsion}
Let $(M,J,g)$ be a $(J^2=\pm 1)$-metric manifold and let $\nabla^{\mathrm{a}}$ be a covariant derivative on $M$ adapted to $(J,g)$.
Then the torsion tensor $\mathrm{T}^{\mathrm{a}}$ of $\nabla^{\mathrm{a}}$  satisfies
\[
N_J(X,Y) = J\mathrm{T}^{\mathrm{a}}(JX,Y)+J\mathrm{T}^{\mathrm{a}}(X,JY)-\alpha\mathrm{T}^{\mathrm{a}}(X,Y)-\mathrm{T}^{\mathrm{a}}(JX,JY),   \quad \forall X,Y \in \mathfrak{X} (M).
\]
\end{lema}

{\bf Proof.} Let $S=\nabla^{\mathrm{a}}-\nabla^{\mathrm{g}}$ be the potential tensor of $\nabla^{\mathrm{a}}$. By Lemma \ref{teor:natural} one has
\begin{eqnarray*}
(\nabla^{\mathrm{g}}_X J) JY = JS(X,JY)-\alpha S(X,Y), &\quad& (\nabla^{\mathrm{g}}_{JX} J) Y = JS(JX,Y)-S(JX,JY), \\
(\nabla^{\mathrm{g}}_Y J) JX = JS(Y,JX)-\alpha S(Y,X), &\quad& (\nabla^{\mathrm{g}}_{JY} J) X = JS(JY,X)-S(JY,JX),   \quad \forall X, Y \in \mathfrak{X} (M),
\end{eqnarray*}
and according to formulas  (\ref{teor:tora-potencial}) and (\ref{eq:nijenhuis}),  one obtains
\[
N_J(X,Y)=  J\mathrm{T}^{\mathrm{a}}(JX,Y)+J\mathrm{T}^{\mathrm{a}}(X,JY)-\alpha\mathrm{T}^{\mathrm{a}}(X,Y)-\mathrm{T}^{\mathrm{a}}(JX,JY),  \quad \forall X, Y \in \mathfrak{X} (M). \ \blacksquare
\]

\begin{lema}
\label{teor:mascondicionessuficientes}
Let $(M,J,g)$ be a $(J^2=\pm 1)$-metric manifold and let $\nabla^{\mathrm{a}}$ be a covariant derivative on $M$ adapted to $(J,g)$. The following relations hold:

\begin{enumerate}
\renewcommand*{\theenumi}{\roman{enumi})}
\renewcommand*{\labelenumi}{\theenumi}

\item If $\mathrm{T}^{\mathrm{a}}(JX,JY)=(-\alpha) \mathrm{T}^{\mathrm{a}}(X,Y)$, $\forall X, Y \in \mathfrak{X} (M)$, then the Nijenhuis tensor of $J$ vanishes.

\item If $J \mathrm{T}^{\mathrm{a}}(JX,Y)=\alpha \mathrm{T}^{\mathrm{a}}(X,Y)$,  $\forall X, Y \in \mathfrak{X} (M)$, then the Nijenhuis tensor of $J$ vanishes.
\end{enumerate}
\end{lema}

{\bf Proof.} Given $X, Y$ vector fields on $M$, if $\mathrm{T}^{\mathrm{a}}(JX,JY)=(-\alpha) \mathrm{T}^{\mathrm{a}}(X,Y)$ then $\mathrm{T}^{\mathrm{a}}(JX,Y)=-\mathrm{T}^{\mathrm{a}}(X,JY)$,
and if  $J \mathrm{T}^{\mathrm{a}}(JX,Y)=\alpha \mathrm{T}^{\mathrm{a}}(X,Y)$ then $J\mathrm{T}^{\mathrm{a}}(X,JY)=\mathrm{T}^{\mathrm{a}}(JX,JY)$.

The result follows from the above equalities and  Lemma  \ref{teor:nijenhuis-torsion}. $\blacksquare$

\subsection{The second Nijenhuis tensor and quasi-K\"{a}hler type manifolds}

As we have shown in Lemma \ref{teor:kahler-caracterizacion}, K\"{a}hler type manifolds are those manifolds for which the Levi Civita connection is natural respect to the $(\alpha ,\varepsilon )$-structure. By formula (\ref{eq:nijenhuis}) we know that the $\alpha$-structure $J$ of a K\"{a}hler type manifold is integrable. We are looking for a new tensor which allows to characterize quasi-K\"{a}hler type manifolds. This tensor will be called the second Nijenhuis tensor. Let us begin introducing the tensor, studying its main properties and, after that, remembering the notion of quasi-K\"{a}hler type manifold and comparing with the vanishing of the second Nijenhuis tensor.

Taking in mind formula (\ref{eq:nijenhuis}) for the Nijenhuis tensor, we introduce the following:

\begin{defin}
\label{teor:njae}
Let $(M,J,g)$ be a  $(J^2=\pm 1)$-metric manifold. The second  Nijenhuis tensor of $(J,g)$ is the tensor field of type $(1,1)$ given by
\[
\widetilde N_J^{\alpha\varepsilon}(X,Y)=(\nabla^{\mathrm{g}}_X J) JY+ \alpha \varepsilon \big( (\nabla^{\mathrm{g}}_{JX} J)Y +(\nabla^{\mathrm{g}}_Y J)JX \big)+(\nabla^{\mathrm{g}}_{JY}J) X,   \quad \forall X, Y \in \mathfrak{X} (M).
\]
\end{defin}

Obviously, this definition is expressed in terms of $\nabla^{\mathrm{g}} J$, which is one of the aims of this section. Observe that the above definition depends on the value $\alpha \varepsilon$ and not just of the $\alpha$-structure. The next properties follow in a direct way.

\begin{lema}
\label{teor:propiedadesNJ0}
Let $(M,J,g)$ be a $(J^2=\pm 1)$-metric manifold. The second Nijenhuis tensor of $(J,g)$ satisfies:

\begin{enumerate}
\renewcommand*{\theenumi}{\roman{enumi})}
\renewcommand*{\labelenumi}{\theenumi}
\item
$\widetilde N_J^{\alpha\varepsilon}(Y,X)=\alpha\varepsilon\widetilde N_J^{\alpha\varepsilon}(X,Y), 
\widetilde N_J^{\alpha\varepsilon}(JX,JY)=\varepsilon \widetilde N_J^{\alpha\varepsilon}(X,Y), 
\widetilde N_J^{\alpha\varepsilon}(JX,Y)=\alpha\varepsilon\widetilde N_J^{\alpha\varepsilon}(X,JY)$, for all vector fields  $X, Y$ on $M$.
\item
\begin{eqnarray*}
g(\widetilde N_J^{\alpha\varepsilon}(X,Y),JZ) &=& -\varepsilon g((\nabla^{\mathrm{g}}_{X} J) Y,Z)   -g ((\nabla^{\mathrm{g}}_{JX} J) JY, Z)-\alpha g((\nabla^{\mathrm{g}}_{Y} J) X, Z) - \alpha \varepsilon g((\nabla^{\mathrm{g}}_{JY} J) JX,Z),
\\
g(\widetilde N_J^{\alpha\varepsilon}(JX,Y),JZ) &=& -\varepsilon g((\nabla^{\mathrm{g}}_{JX} J) Y,Z)   -\alpha( g ((\nabla^{\mathrm{g}}_{X} J) JY, Z)+ g((\nabla^{\mathrm{g}}_{Y} J) JX, Z)) - \varepsilon g((\nabla^{\mathrm{g}}_{JY} J) X,Z),
\end{eqnarray*}
for all vector fields $X, Y, Z$ on $M$.
\end{enumerate}
\end{lema}

Definition of quasi-K\"{a}hler type manifold depends on the geometry we are considering. The four geometries of $(\alpha ,\varepsilon )$-structures have had each own development. In order to have a common presentation of the notion we must distinguish the cases $\alpha\varepsilon=1$ and  $\alpha\varepsilon=-1$. Moreover, this will be useful to compare quasi-K\"{a}hler type manifolds with the vanishing of  the second Nijenhuis tensor.

\begin{defin}
\label{teor:quasikahler+1}
Let $(M,J,g)$ be a $(J^2=\pm1)$-metric manifold with $\alpha\varepsilon=1$. It is said to be a  quasi-K\"{a}hler type manifold if
\[
g((\nabla^{\mathrm{g}}_X J)Y, Z)+g((\nabla^{\mathrm{g}}_Y J)Z, X)+g((\nabla^{\mathrm{g}}_Z J)X,Y)=0,   \quad \forall X, Y, Z\in \mathfrak{X} (M).
\]
\end{defin}

According to Definition \ref{teor:njae} the second Nijenhuis tensor for a $(J^2=\pm1)$-metric manifold with $\alpha\varepsilon=1$ is
\begin{equation}
\widetilde N_J^{1}(X,Y)=(\nabla^{\mathrm{g}}_X J) JY +(\nabla^{\mathrm{g}}_{JX} J)Y +(\nabla^{\mathrm{g}}_Y J)JX+(\nabla^{\mathrm{g}}_{JY}J) X,   \quad \forall X, Y \in \mathfrak{X} (M).
\label{eq:njt}
\end{equation}

The characterization of  quasi-K\"{a}hler type manifolds in terms of the second Nijenhuis tensor field has been obtained for each of the two geometries:

\begin{prop}[{\cite{ganchev-borisov}, \cite{staikova}}]
\label{teor:caracterizacion-quasikahler+1}
Let $(M,J,g)$ be a $(J^2=\pm 1)$-metric manifold with $\alpha \varepsilon =1$. The following conditions are equivalent:

\begin{enumerate}
\renewcommand*{\theenumi}{\roman{enumi})}
\renewcommand*{\labelenumi}{\theenumi}
\item The second Nijenhuis tensor  of $(J,g)$ vanishes.

\item The manifold $(M,J,g)$ is a  quasi-K\"{a}hler type manifold.
\end{enumerate}
\end{prop}

Quasi-K\"{a}hler type manifolds with $\alpha \varepsilon=1$ have been studied in several papers as  \cite{mekerov} and \cite{teofilova}, in the case of $(1,1)$-structures, and  $(-1,-1)$-structures, respectively. These manifolds correspond to the class $\mathcal W_3$ in the classification of $(\alpha,\varepsilon )$-structures with $\alpha \varepsilon =1$, and it is the unique class of the basic ones characterized by the non-integrability of the $\alpha$-structure  $J$ (see \cite{ganchev-borisov} and \cite{staikova}, where  almost Norden and almost-product Riemannian manifolds with null trace are classified, and Proposition \ref{teor:quasikahler-J-integrable}).
\bigskip

In the case $\alpha\varepsilon=-1$ expression of the second Nijenhuis tensor of $(J,g)$ is
\begin{equation}
\widetilde N_J^{-1}(X,Y)=(\nabla^{\mathrm{g}}_X J) JY-(\nabla^{\mathrm{g}}_{JX} J)Y -(\nabla^{\mathrm{g}}_Y J)JX+(\nabla^{\mathrm{g}}_{JY}J) X,   \quad \forall X, Y \in \mathfrak{X} (M),
\label{eq:njd}
\end{equation}
and the definition of quasi-K\"ahler type manifolds is:
\begin{defin}
Let $(M,J,g)$ be a $(J^2=\pm1)$-metric manifold with $\alpha\varepsilon=-1$. It is said to be a  quasi-K\"{a}hler type manifold if
\begin{equation}
(\nabla^{\mathrm{g}}_X J)JY-(\nabla^{\mathrm{g}}_{JX} J) Y=0,   \quad \forall X, Y\in \mathfrak{X} (M).
\label{eq:quasikahler-1}
\end{equation}
\end{defin}

Then we obtain the following characterization of quasi-K\"{a}hler type manifolds by means of the second Nijenhuis tensor:
\begin{prop}
\label{teor:caracterizacion-quasikahler-1}
Let $(M,J,g)$ be a $(J^2=\pm 1)$-metric manifold with $\alpha \varepsilon =-1$. The following conditions are equivalent:
\begin{enumerate}
\renewcommand*{\theenumi}{\roman{enumi})}
\renewcommand*{\labelenumi}{\theenumi}
\item The second Nijenhuis tensor of $(M,J,g)$ vanishes.

\item $(\nabla^{\mathrm{g}}_X \Phi) (Y,Z)-\alpha (\nabla^{\mathrm{g}}_{JX} \Phi) (JY,Z)=0$, for all vector fields  $X, Y, Z$ on  $M$.

\item $(\nabla^{\mathrm{g}}_X J) Y -\alpha (\nabla^{\mathrm{g}}_{JX} J) JY=0$, for all vector fields $X,Y$ on $M$.

\item The manifold $(M,J,g)$ is a  quasi-K\"{a}hler type manifold.

\item $(\nabla^{\mathrm{g}}_X J) X - \alpha (\nabla^{\mathrm{g}}_{JX} J) JX =0$, for all vector field $X$ on $M$. 
\end{enumerate}
\end{prop}

{\bf Proof.}

$i)\Rightarrow ii)$ As $\alpha\varepsilon=-1$, according to Lemmas \ref{teor:tensorNJ}  and \ref{teor:propiedadesNJ0}, one has
\[
g( \widetilde N_J^{-1} (X,Y),JZ)-g(\widetilde N_J^{-1}(Y,Z),JX)+g(\widetilde N_J^{-1}(Z,X),JY)=  2(\alpha g((\nabla^{\mathrm{g}}_X J) Y, Z) -  g((\nabla^{\mathrm{g}}_{JX} J) JY,Z)),
\]
for all vector fields $X,Y,Z$ on $M$. Then, the vanishing of the second Nijenhuis tensor of $(J,g)$ implies
\[
(\nabla^{\mathrm{g}}_X \Phi) (Y,Z)-\alpha (\nabla^{\mathrm{g}}_{JX} \Phi) (JY,Z)=0,   \quad \forall X, Y, Z \in \mathfrak{X} (M).
\]

$ii)\Rightarrow iii)$ It follows in a direct way from formula (\ref{eq:pp-tensorfundamental}).

$iii) \Rightarrow iv)$ Evaluate the expression in  $(X,JY)$.

$iv) \Rightarrow v)$ As $(M,J,g)$ is a quasi-K\"{a}hler type manifold, given $X, Y$ vector fields on $M$ such that $Y=JX$ from (\ref{eq:quasikahler-1}) one obtains
\[
\alpha (\nabla^{\mathrm{g}}_X J) X -  (\nabla^{\mathrm{g}}_{JX} J) JX =0 \Leftrightarrow (\nabla^{\mathrm{g}}_X J) X - \alpha (\nabla^{\mathrm{g}}_{JX} J) JX =0.
\]

$v)\Rightarrow i)$ Given $X, Y$ vector fields on $M$ one has 
\begin{eqnarray*}
(\nabla^{\mathrm{g}}_{X+JY} J) (X+JY) &=& (\nabla^{\mathrm{g}}_X J) X + (\nabla^{\mathrm{g}}_X J) JY + (\nabla^{\mathrm{g}}_{JY} J) X + (\nabla^{\mathrm{g}}_{JY} J)JY, \\
\alpha (\nabla^{\mathrm{g}}_{J(X+JY)} J) J(X+JY) &=&  \alpha (\nabla^{\mathrm{g}}_{JX + \alpha Y} J) (JX+\alpha Y)= \alpha (\nabla^{\mathrm{g}}_{JX} J) JX+(\nabla^{\mathrm{g}}_{JX} J) Y +  (\nabla^{\mathrm{g}}_Y J) JX + \alpha (\nabla^{\mathrm{g}}_{Y} J)Y.
\end{eqnarray*}
Subtracting both equations and taking into account formula (\ref{eq:njd}) one obtains the result. $\blacksquare$
\bigskip

Manifolds having $\alpha \varepsilon=-1$ correspond to almost Hermitian, i.e., $(-1,1)$, manifolds and almost para-Hermitian, i.e., $(1,-1)$, manifolds. In the case of almost Hermitian, quasi-K\"ahler manifolds are introduced in \cite{gray-hervella} as manifolds satisfying condition $ii)$ of the above result, with $\alpha=-1$. In the almost para-Hermitian case quasi-K\"ahler type manifolds are introduced in \cite{chursinetal} as we have written in the corresponding definition. A characterization of quasi-K\"ahler type manifolds in almost Hermitian and almost para-Hermitian geometries in terms of the second Nijenhuis tensor has not been previously obtained. In the other two geometries the corresponding characterizations had been obtained, as we have indicated in Proposition \ref{teor:caracterizacion-quasikahler+1}, in terms of the second Nijenhuis tensor which  is given by formula (\ref{eq:njt}), in this case of $\alpha \varepsilon =1$. In order to have Proposition \ref{teor:caracterizacion-quasikahler-1} we have had to obtain an expression of the second Nijenhuis tensor for the case $\alpha \varepsilon =-1$ compatible with that for the case $\alpha \varepsilon =1$. This was done in Definition \ref{teor:njae}.
\bigskip

One obtains the following technical result: 
\begin{lema}
Let $(M,J,g)$ be a $(J^2=\pm 1)$-metric manifold with $\alpha \varepsilon =-1$. The following conditions are equivalent:
\begin{enumerate}
\renewcommand*{\theenumi}{\roman{enumi})}
\renewcommand*{\labelenumi}{\theenumi}
\item $(\nabla^{\mathrm{g}}_X J) JY-(\nabla^{\mathrm{g}}_{JX} J)Y +(\nabla^{\mathrm{g}}_Y J)JX-(\nabla^{\mathrm{g}}_{JY}J) X=0$, for all vector fields $X, Y$ on $M$.

\item $(\nabla^{\mathrm{g}}_X J) JX - (\nabla^{\mathrm{g}}_{JX} J) X =0$, for all vector field $X$ on $M$. 

\item $(\nabla^{\mathrm{g}}_X J) X - \alpha (\nabla^{\mathrm{g}}_{JX} J) JX =0$, for all vector field $X$ on $M$. 
\end{enumerate}
\end{lema}

{\bf Proof.}

$i) \Rightarrow ii)$ Given $X, Y$ vector fields on $M$ such that $X=Y$ one obtains $(\nabla^{\mathrm{g}}_X J)JX - (\nabla^{\mathrm{g}}_{JX} J) X=0$.

$ii) \Rightarrow iii)$ Given $X, Y$ vector fields on $M$ one has
\[
g((\nabla^{\mathrm{g}}_X J)JX, JY) - g((\nabla^{\mathrm{g}}_{JX} J) X, JY)=0, 
\]
then taking into account (\ref{eq:tensorNJ3}) one obtains
\[
\alpha g((\nabla^{\mathrm{g}}_X J)X,Y) - g((\nabla^{\mathrm{g}}_{JX} J) JX,Y)=0, 
\]
thus 
\[
\alpha (\nabla^{\mathrm{g}}_X J)X - (\nabla^{\mathrm{g}}_{JX} J) JX=0 \Leftrightarrow (\nabla^{\mathrm{g}}_X J)X -\alpha  (\nabla^{\mathrm{g}}_{JX} J) JX=0.
\]

$iii) \Rightarrow i)$ In these conditions the manifold $(M,J,g)$ is a quasi-K\"{a}hler manifold (see Proposition \ref{teor:caracterizacion-quasikahler-1} $v)$), then  one has
\[
(\nabla^{\mathrm{g}}_X J) JY-(\nabla^{\mathrm{g}}_{JX} J)Y=0, \quad (\nabla^{\mathrm{g}}_Y J)JX-(\nabla^{\mathrm{g}}_{JY}J) X=0, \quad \forall X, Y \in \mathfrak{X} (M),
\] 
adding both equations one obtains the result. $\blacksquare$

\begin{obs}
\label{teor:ttnj4}
Then the following condition
\[
(\nabla^{\mathrm{g}}_X J) JY-(\nabla^{\mathrm{g}}_{JX} J)Y +(\nabla^{\mathrm{g}}_Y J)JX-(\nabla^{\mathrm{g}}_{JY}J) X=0, \quad \forall X, Y \in \mathfrak{X} (M),
\]
also characterizes the $(J^2=\pm1)$-metric manifolds of quasi-K\"{a}hler type in the case $\alpha \varepsilon =-1$.
\end{obs}

Now we prove a result relating integrable $J$-structures, quasi-K\"ahler type and K\"ahler type manifolds:

\begin{prop}
\label{teor:quasikahler-J-integrable}
Let $(M,J,g)$ be a $(J^2=\pm1)$-metric manifold. If $(M,J,g)$ is a quasi-K\"{a}hler type manifold such that its $\alpha$-structure $J$ is integrable then $(M,J,g)$ is a K\"{a}hler type manifold.
\end{prop}

{\bf Proof.} In order to prove the result we must distinguish the cases $\alpha \varepsilon =\pm 1$.

 Assuming $\alpha\varepsilon=-1$,  and according to  Lemma \ref{teor:caracterizacionintegrabilidad-1}, the $\alpha$-structure $J$ is integrable if and only if
\[
(\nabla^{\mathrm{g}}_X J) Y +\alpha (\nabla^{\mathrm{g}}_{JX} J) JY =0,   \quad \forall X, Y \in \mathfrak{X} (M).
\]
As $(M,J,g)$ is a quasi-K\"{a}hler type manifold, and according to Proposition \ref{teor:caracterizacion-quasikahler-1} one has
\[
(\nabla^{\mathrm{g}}_X J) Y -\alpha (\nabla^{\mathrm{g}}_{JX} J) JY =0,   \quad \forall X, Y \in \mathfrak{X} (M),
\]
thus proving  $\nabla^{\mathrm{g}} J=0$, adding both equalities.

In the case $\alpha\varepsilon=1$, the vanishing of both the Nijenhuis tensor and the second Nijenhuis tensor implies
\[
(\nabla^{\mathrm{g}}_X J) JY + (\nabla^{\mathrm{g}}_{JX} J) Y =0,   \quad \forall X, Y \in \mathfrak{X} (M).
\]
Adding the expression in Definition \ref{teor:quasikahler+1} valued  in $(X,Y,JZ)$ and the expression of the second property of Lemma \ref{teor:caracterizacionintegrabilidad+1} one has
\[
2 g((\nabla^{\mathrm{g}}_X J)Y, JZ) + g((\nabla^{\mathrm{g}}_Z J)JY,X)+g((\nabla^{\mathrm{g}}_{JZ} J)Y,X)=0,   \quad \forall X, Y, Z \in \mathfrak{X} (M),
\]
which implies $\nabla^{\mathrm{g}} J=0$ taking in mind the above equality.

In both cases $\alpha \varepsilon = \pm 1$ we have proved $\nabla^{\mathrm{g}} J=0$. As $\nabla^{\mathrm{g}} g=0$, then $\nabla^{\mathrm{g}}$ is adapted to  $(J,g)$ and, according to Lemma \ref{teor:kahler-caracterizacion},  $(M,J,g)$ is a  K\"{a}hler type manifold. $\blacksquare$
\bigskip

We finish this section recalling briefly the nearly K\"{a}hler type manifolds in the case of $\alpha\varepsilon=-1$. They were introduced by Gray in the almost Hermitian case (see \cite{gray}) and correspond to a class in the classification of almost Hermitian manifolds of Gray and Hervella (see \cite{gray-hervella}). In the almost para-Hermitian case the analogous class also appears in the classification of Gadea and Mu\~{n}oz Masqu\'{e} (see \cite{gadea}), where two of the eight classes are the so-called $(+)$-nearly para-K\"{a}hlerian and $(-)$-nearly para-K\"{a}hlerian manifolds. 

\begin{defin}[{\cite{gray}, \cite{ivanov}}]
Let $(M,J,g)$ be a $(J^2=\pm1)$-metric manifold with $\alpha \varepsilon =-1$. Then $(M,J,g)$ is said to be a  nearly K\"{a}hler type manifold if the following relation holds
\[
(\nabla^{\mathrm{g}}_X J) X=0, \quad \forall X \in \mathfrak{X} (M).
\]
\end{defin}

\begin{lema}
\label{teor:nearlykahler}
Let $(M,J,g)$ be a $(J^2=\pm1)$-metric manifold such that $\alpha \varepsilon =-1$. The following conditions are equivalent:
\begin{enumerate}
\renewcommand*{\theenumi}{\roman{enumi})}
\renewcommand*{\labelenumi}{\theenumi}

\item $(\nabla^{\mathrm{g}}_X J) X=0$, for all vector field $X$ on $M$.

\item $(\nabla^{\mathrm{g}}_X J) Y + (\nabla^{\mathrm{g}}_Y J)X=0$, for all vector fields $X, Y$ on $M$.
\end{enumerate}
\end{lema}

{\bf Proof.} Trivial.

\begin{prop}
\label{teor:nearly-quasi}
Let $(M,J,g)$ be a $(J^2=\pm1)$-metric manifold with $\alpha \varepsilon =-1$. If $(M,J,g)$ is a nearly K\"{a}hler type manifold then is also a quasi-K\"{a}hler type manifold and the Nijenhuis tensor of $J$ satisfies the following relation
\[
N_J(X,Y)= 2((\nabla^{\mathrm{g}}_X J)JY + (\nabla^{\mathrm{g}}_{JX} J) X), \quad \forall X, Y \in \mathfrak{X} (M).
\]
\end{prop}

{\bf Proof.} As $(M,J,g)$ is a nearly K\"{a}hler manifold, given $X, Y$ vector fields on $M$ by the previous Lemma one has
\begin{equation}
(\nabla^{\mathrm{g}}_X J) JY + (\nabla^{\mathrm{g}}_{JY} J) X=0, \quad (\nabla^{\mathrm{g}}_{JX} J) Y + (\nabla^{\mathrm{g}}_{Y} J) JX=0. 
\label{eq:nearlykahler2}
\end{equation}
Substracting both equations one obtains
\[
(\nabla^{\mathrm{g}}_X J) JY - (\nabla^{\mathrm{g}}_{JX} J) Y - (\nabla^{\mathrm{g}}_{Y} J) JX + (\nabla^{\mathrm{g}}_{JY} J) X=0, \quad \forall X, Y \in \mathfrak{X} (M), 
\]
then according to (\ref{eq:njt}) and Proposition \ref{teor:caracterizacion-quasikahler-1} $(M,J,g)$ is a quasi-K\"{a}hler type manifold.  

Taking into account the equalities (\ref{eq:nearlykahler2}) one also obtains
\begin{eqnarray*}
N_J(X,Y) &=&(\nabla^{\mathrm{g}}_X J) JY + (\nabla^{\mathrm{g}}_{JX} J) Y - (\nabla^{\mathrm{g}}_Y J) JX -(\nabla^{\mathrm{g}}_{JY} J) X\\
               &=& (\nabla^{\mathrm{g}}_X J) JY + (\nabla^{\mathrm{g}}_{JX} J) Y + (\nabla^{\mathrm{g}}_{JX} J) Y +(\nabla^{\mathrm{g}}_{X} J) JY\\
               &=& 2 ((\nabla^{\mathrm{g}}_X J) JY + (\nabla^{\mathrm{g}}_{JX} J) Y), \quad \forall X, Y \in \mathfrak{X} (M). \blacksquare
\end{eqnarray*}

\section{Distinguished connections on $(J^2=\pm1)$-metric manifolds}

This section and the following one are the core of the paper. We consider a $(J^2=\pm1)$-metric manifold. Our main aims are:
\begin{enumerate}
\item Studying the distinguished  connections, namely, the first canonical, the Chern, the well adapted, and the Kobayashi-Nomizu and the Yano connections, and connections with skew-symmetric torsion tensor when they can be defined.
\item Characterizing the above connections by the vanishing of suitable tensor fields defined in the previous section, when possible, and obtaining properties of the torsion tensor.

\item Characterizing the coincidence among connections, when possible.
\end{enumerate}

We present these connections in the quoted order and we study simultaneously the characterization properties. In particular we will prove:

\begin{itemize}
\item The first canonical connection $\nabla ^{0}$ can be defined in any $(J^2=\pm1)$-metric manifold and it is adapted to the structure.

\item The Chern connection $\nabla^{\mathrm{c}}$ can be defined in  $(J^2=\pm1)$-metric manifold with $\alpha \varepsilon =-1$ and it is adapted to the structure. There is no a definition for the case $\alpha \varepsilon =1$. Assuming $\alpha \varepsilon =-1$,  it will be proved that  $\nabla ^{0}=\nabla^{\mathrm{c}}$ if and only if the manifold is quasi-K\"{a}hler.

\item The well adapted connection $\nabla^{\mathrm{w}}$ can be defined in any $(J^2=\pm1)$-metric manifold and it is adapted to the structure. A characterization of $\nabla ^{0}=\nabla^{\mathrm{w}}$ will be obtained.

\item The Kobayashi-Nomizu connection $\nabla^{\mathrm{kn}}$ can be defined in any $(J^2=\pm1)$-metric manifold but it is not natural in general. When a  Kobayashi-Nomizu  connection is natural will be completely characterized.

\item The Yano connection $\nabla^{\mathrm{y}}$    can be defined in any $(J^2=\pm1)$-metric manifold but it is not natural in general (let us remember that it is not even adapted to $(M,J)$ in the general case). We will characterize the case when a Yano connection is adapted.

\item Connections with skew-symmetric torsion tensor $\nabla^{\mathrm{sk}}$ can be defined in any $(J^2=\pm1)$-metric manifold but they are  not natural in general. Adapted connections will be characterized. These connections are not uniquely defined.
\end{itemize}

Besides, one can consider the Levi Civita connection $\nabla^{\mathrm{g}}$, which is natural if and only if the manifold is a  K\"{a}hler type manifold (Lemma \ref{teor:kahler-caracterizacion}). In fact, $\nabla ^{0}$, $\nabla^{\mathrm{kn}}$ and $\nabla^{\mathrm{y}}$ are connections uniquely defined from the Levi Civita connection, while the other $\nabla^{\mathrm{c}}$, $\nabla^{\mathrm{w}}$ and $\nabla^{\mathrm{sk}}$ are connections defined imposing being adapted and satisfying some conditions on the torsion tensor. In this second case one should prove uniqueness, if there exists.

In the last section of the paper we will define a 1-parameter family of adapted connections. This family will contain other distinguished connections such as the Bismut connection $\nabla^{\mathrm{b}}$, as we will show.

\subsection{The first canonical connection}

The covariant derivative $\nabla ^{0}$ of this connection was introduced in  Definition \ref{teor:first-connection} (in the almost Hermitian case the definition was given in classical and seminal papers as  \cite{lichnerowicz2} and  \cite{gauduchon}). According to formula   (\ref{teor:tora-potencial}), its torsion tensor has the following expression
\begin{equation}
\mathrm{T}^{0}(X,Y)=\frac{(-\alpha)}{2} ((\nabla^{\mathrm{g}}_X J)JY-(\nabla^{\mathrm{g}}_Y J) JX),   \quad \forall X,Y \in \mathfrak{X} (M).
\label{eq:torsion0}
\end{equation}

As $\nabla ^{0}$ is a natural connection,  Lemma \ref{teor:nijenhuis-torsion} is valid for it, thus establishing a link between the torsion tensor $\mathrm{T}^{0}$ and  the Nijenhuis tensor. The following results show other properties of $\mathrm{T}^{0}$.

\begin{prop}
Let $(M,J,g)$ be a $(J^2=\pm1)$-metric manifold. The following relation holds:
\begin{eqnarray}
\mathrm{T}^{0}(JX,JY) +\alpha \mathrm{T}^{0} (X,Y) &=& -\frac{1}{2} N_J(X,Y), \label{eq:toro+alpha}\\
\mathrm{T}^{0}(JX,JY) -\alpha \mathrm{T}^{0} (X,Y) &=& \frac{1}{2} ((\nabla^{\mathrm{g}}_X J)JY- (\nabla^{\mathrm{g}}_{JX} J)Y -(\nabla^{\mathrm{g}}_Y J)JX)+(\nabla^{\mathrm{g}}_{JY} J)X), \label{eq:toro-alpha}
\end{eqnarray}
for all vector fields $X,Y$ on $M$.
\end{prop}

{\bf Proof. }  According to formulas (\ref{eq:torsion0}) and (\ref{eq:nijenhuis}) we obtain for all vector fields $X,Y$ on $M$
\begin{eqnarray*}
\mathrm{T}^{0}(JX,JY) +\alpha \mathrm{T}^{0} (X,Y) &=& -\frac{1}{2} ((\nabla^{\mathrm{g}}_X J)JY+ (\nabla^{\mathrm{g}}_{JX} J)Y -(\nabla^{\mathrm{g}}_Y J)JX)-(\nabla^{\mathrm{g}}_{JY} J)X) =-\frac{1}{2} N_J(X,Y),\\
\mathrm{T}^{0}(JX,JY) -\alpha \mathrm{T}^{0} (X,Y) &=& \frac{1}{2} ((\nabla^{\mathrm{g}}_X J)JY- (\nabla^{\mathrm{g}}_{JX} J)Y -(\nabla^{\mathrm{g}}_Y J)JX)+(\nabla^{\mathrm{g}}_{JY} J)X),
\end{eqnarray*}
thus proving the result. $\blacksquare$
\bigskip

Observe that in the case $\alpha\varepsilon=-1$, taking into account formula (\ref{eq:njd}),  one has $
\mathrm{T}^{0}(JX,JY) -\alpha \mathrm{T}^{0} (X,Y) =\frac{1}{2} \widetilde N_J^{-1}(X,Y)$, for all vector fields  $X, Y$ on $M$. Then one easily proves:

\begin{corol}
\label{teor:toro-integrability}
Let $(M,J,g)$ be a $(J^2=\pm1)$-metric manifold.
\begin{enumerate}
\renewcommand*{\theenumi}{\roman{enumi})}
\renewcommand*{\labelenumi}{\theenumi}

\item If $\alpha=-1$ then $J$ is integrable if and only if $\mathrm{T}^{0}(JX,JY)=\mathrm{T}^{0}(X,Y)$, for all vector fields  $X, Y$ on $M$.

\item If  $\alpha=1$ then $J$  is integrable if and only if $\mathrm{T}^{0}(JX,JY)=-\mathrm{T}^{0}(X,Y)$, for all vector fields  $X, Y$ on $M$.

\item If $\alpha\varepsilon =-1$ then $(M,J,g)$ is a quasi-K\"{a}hler type manifold if and only if $\mathrm{T}^{0}(JX,JY)=\alpha \mathrm{T}^{0}(X,Y)$, for all vector fields  $X, Y$ on $M$.
\end{enumerate}
\end{corol}

\subsection{The Chern connection}

The Chern connection was firstly introduced in the case of almost Hermitian manifolds. In  \cite{racsam} we have extended the connection to the almost para-Hermitian case, recovering the connection defined by Cruceanu and one of us in \cite{etayo}. The following results establish the existence and uniqueness of the Chern connection in a $(J^2=\pm1)$-metric manifold with $\alpha\varepsilon=-1$.

\begin{teor}[{\cite[Theor. 6.3]{racsam}}]
\label{teor:chern-connection}
Let  $(M,J,g)$ be a $(J^2=\pm1)$-metric manifold with  $\alpha\varepsilon =-1$. Then there exists a unique linear connection $\nabla^{\mathrm{c}}$ in $M$ reducible to the $G_{(\alpha,\varepsilon)}$-structure defined by $(J,g)$  whose torsion tensor  $\mathrm{T}^{c}$ satisfies the following condition
\[
\mathrm{T}^{c}(JX,JY)= \alpha \mathrm{T}^{c}(X,Y),   \quad \forall X, Y \in \mathfrak{X} (M).
\]
This connection is called the Chern connection of $(M,J,g)$.
\end{teor}

According to the second condition in Lemma \ref{teor:mascondicionessuficientes} to the covariant derivative $\nabla^{\mathrm{c}}$ of the Chern connection, the condition
\[
J\mathrm{T}^{c}(JX,Y)= \alpha\mathrm{T}^{c}(X,Y),   \quad \forall X,Y \in \mathfrak{X} (M),
\]
determines the Chern connection in the case $J$ is integrable. This property has been taken in  \cite[Theor. 3.5]{ivanov} in order to introduce the Chern connection in para-Hermitian manifolds. Our point of view is more general, because it also includes the non-integrable case.

\bigskip

The following result characterizes the identity $\nabla^{0}=\nabla^{\mathrm{c}}$:
\begin{prop}
\label{teor:n0=nc}
Let $(M,J,g)$ be a  $(J^2=\pm1)$-metric manifold with $\alpha\varepsilon=-1$. Then the first canonical connection and the   Chern connection of $(M,J,g)$ coincide if and only if the second   Nijenhuis tensor of $(M,J,g)$ vanishes.
\end{prop}

{\bf Proof.}  By formula (\ref{eq:toro-alpha})  one has
\[
\mathrm{T}^{0}(JX,JY)= \alpha \mathrm{T}^{0}(X,Y) \Leftrightarrow \widetilde N_J^{-1}(X,Y)=0,   \quad \forall X, Y \in \mathfrak{X} (M). \ \blacksquare
\]

According to Proposition \ref{teor:caracterizacion-quasikahler-1} we can also have written $\nabla^{0}=\nabla^{\mathrm{c}}$ if and only if $(M,J,g)$ is a quasi-K\"{a}hler type manifold.

\subsection{The well adapted connection}

In \cite{racsam} we have deeply studied this connection. It is an adapted connection to $(M,J,g)$. It is also a functorial connection and it is the most natural connection in the following sense: the $G_{(\alpha,\varepsilon )}$-structure defined by $(M,J,g)$ is integrable if and only if the torsion and the curvature tensors of the well adapted connection vanish.

 As in the case of the Chern connection, the well adapted connection can be defined as the unique connection satisfying a condition about its torsion tensor:

\begin{teor}[{\cite[Theor. 4.4]{racsam}}]
 \label{teor:bienadaptada-ae-estructura}
Let  $(M,J,g)$ be a $(J^2=\pm1)$-metric manifold. Then there exists a unique linear connection $\nabla^{\mathrm{w}}$ in $M$ reducible to the $G_{(\alpha,\varepsilon)}$-structure defined by $(J,g)$  whose torsion tensor  $\mathrm{T}^{\mathrm{w}}$ satisfies the following condition
 \begin{equation}
g(\mathrm{T}^{\mathrm{w}}(X,Y),Z)-g(\mathrm{T}^{\mathrm{w}}(Z,Y),X) = -\varepsilon (g(\mathrm{T}^{\mathrm{w}}(JX,Y),JZ)-g(\mathrm{T}^{\mathrm{w}}(JZ,Y),JX)),  \quad \forall X,Y, Z \in \mathfrak{X} (M).
\label{eq:welladapted}
\end{equation}
This connection is called the well adapted connection of $(M,J,g)$.
\end{teor}

It is known (see \cite[Theor. 5.2]{racsam}) that the well adapted connection and the Levi Civita connection coincide if and only if $(M,J,g)$ is a K\"{a}hler type manifold, or, equivalently if $\nabla ^{\mathrm{g}}$ is a natural covariant derivative (Lemma \ref{teor:kahler-caracterizacion}).
\bigskip

We establish three results about the relation of the well adapted connection with $\nabla ^{0}$, with the integrability of $J$ and with the K\"{a}hler condition.

\begin{teor}
\label{teor:nabla0=nablaw}
Let $(M,J,g)$ be a $(J^2=\pm1)$-metric manifold. Then the first canonical connection and the well adapted connection  coincide if and only if

\begin{enumerate}
\renewcommand*{\theenumi}{\roman{enumi})}
\renewcommand*{\labelenumi}{\theenumi}

\item $(M,J,g)$ quasi-K\"{a}hler type manifold, in the case  $\alpha\varepsilon=-1$.

\item The $\alpha$-structure $J$ is integrable, in the case  $\alpha\varepsilon=1$.
\end{enumerate}

\end{teor}

{\bf Proof.} For an adapted covariant derivative $\nabla^{\mathrm{a}}$ with torsion tensor $\mathrm{T}^{\mathrm{a}}$, let us consider the tensor field ${\mathcal F} (\nabla^{\mathrm{a}})$ of type (0,3) defined as
\[
{\mathcal F} (\nabla^{\mathrm{a}}, X, Y, Z) = g(\mathrm{T}^{\mathrm{a}}(X,Y),Z)- g(\mathrm{T}^{\mathrm{a}}(Z,Y),X) + \varepsilon (g(\mathrm{T}^{\mathrm{a}}(JX,Y),JZ)- g(\mathrm{T}^{\mathrm{a}}(JZ,Y),JX)),
\]
for all vector fields $X,Y,Z$ on $M$. By Theorem \ref{teor:bienadaptada-ae-estructura} we know that ${\mathcal F} (\nabla^{\mathrm{a}})$ vanishes if and only if $\nabla^{\mathrm{a}} =\nabla ^{\mathrm{w}}$. Thus, we want to calculate the tensor field ${\mathcal F} (\nabla ^{0})$ corresponding to the first canonical connection $\nabla ^{0}$ in order to characterize ${\mathcal F} (\nabla ^{0})=0$.

The following relations hold:
\begin{enumerate}
\renewcommand*{\theenumi}{\roman{enumi})}
\renewcommand*{\labelenumi}{\theenumi}

\item If $\alpha\varepsilon =-1$ then
\begin{equation}
{\mathcal F} (\nabla^{0}, X, Y, Z) = \frac{\alpha}{2} g(\widetilde N_J^{-1}(X,Z),Y),   \quad \forall X, Y, Z \in \mathfrak{X} (M).
\label{eq:functorial-primeracanonica-1}
\end{equation}

\item If  $\alpha \varepsilon =1$ then
\begin{equation}
{\mathcal F} (\nabla^{0}, X, Y, Z) = \frac{\alpha}{2} g(N_J(X,Z),Y),    \quad \forall X, Y, Z \in \mathfrak{X} (M).
\label{eq:functorial-primeracanonica+1}
\end{equation}
\end{enumerate}

In order to prove the above relations, consider Lemma \ref{teor:tensorNJ} and formula (\ref{eq:torsion0}) which allows to obtain
\[
{\mathcal F} (\nabla^{0}, X, Y, Z) = \frac{\alpha}{2} (g((\nabla^{\mathrm{g}}_X J)JZ,Y)-g((\nabla^{\mathrm{g}}_Z J)JX),Y)) + \frac{\varepsilon}{2} (g((\nabla^{\mathrm{g}}_{JX} J)Z,Y)-g((\nabla^{\mathrm{g}}_{JZ} J)X),Y)),
\]
for all vector fields $X, Y, Z$ on $M$. Then, taking into account formulas (\ref{eq:njd})  and  (\ref{eq:nijenhuis}) one easily obtains formulas (\ref{eq:functorial-primeracanonica-1}) and (\ref{eq:functorial-primeracanonica+1}).

Finally, taking into account formula (\ref{eq:functorial-primeracanonica-1}) and Proposition \ref{teor:caracterizacion-quasikahler-1}  we obtain the first statement of the present Theorem, and taking into account formula (\ref{eq:functorial-primeracanonica+1}) and the well known fact that the vanishing of the Nijenhuis tensor is equivalent to the integrability of $J$ one has the second statement.
$\blacksquare$
\bigskip

Concerning the integrability of $J$ one obtains:

\begin{prop} Let $(M,J,g)$ be a $(J^2=\pm1)$-metric manifold with $\alpha\varepsilon =1$.  Then the $\alpha$-structure $J$ is integrable if and only if the first canonical connection and the well adapted connection coincide. Besides, $J$ is integrable if and only if
\[
\mathrm{T}^{\mathrm{w}}(JX,JY)=(-\alpha)\mathrm{T}^{\mathrm{w}}(X,Y),   \quad \forall X,Y \in \mathfrak{X} (M).
\]
\end{prop}

{\bf Proof.}  It is a direct consequence of formulas (\ref{eq:functorial-primeracanonica+1}) and (\ref{eq:toro+alpha}). $\blacksquare$
\bigskip

The study of K\"{a}hler condition will be divided in two cases: we will obtain a specific result in the case $\alpha \varepsilon =-1$ and a general result for any $(J^{2}=\pm 1)$-metric manifold.

\begin{prop}
Let $(M,J,g)$ be a  $J^2=\pm 1$-metric manifold with $\alpha\varepsilon =-1$.
\begin{enumerate}
\renewcommand*{\theenumi}{\roman{enumi})}
\renewcommand*{\labelenumi}{\theenumi}

\item If the first canonical connection and the well adapted connection coincide and the $\alpha$-structure $J$ is integrable then $(M,J,g)$ is a K\"{a}hler type manifold.

\item $(M,J,g)$ is a K\"{a}hler type manifold if and only if
\[
\mathrm{T}^{\mathrm{w}}(JX,JY)=(-\alpha)\mathrm{T}^{\mathrm{w}}(X,Y),   \quad \forall X,Y \in \mathfrak{X} (M).
\]
\end{enumerate}
\end{prop}

{\bf Proof.}

$i)$ It follows from Proposition \ref{teor:quasikahler-J-integrable}   and Theorem \ref{teor:nabla0=nablaw}.

$ii)$ We prove both implications.

$\Rightarrow)$ If $(M,J,g)$ is a K\"{a}hler type manifold then $\nabla^{\mathrm{w}}=\nabla^{\mathrm{g}}$, and then $\nabla^{\mathrm{w}}$ is a torsion-free derivative thus obviously satisfying the condition.

$\Leftarrow)$ As $\nabla^{\mathrm{w}}$ is a natural connection, by Lemma  \ref{teor:mascondicionessuficientes} $i)$ one knows $J$ in integrable. Then, by Lemma \ref{teor:caracterizacionintegrabilidad-1} one has

\begin{equation}
(\nabla^{\mathrm{g}}_X \Phi) (Y,Z)= (-\alpha) (\nabla^{\mathrm{g}}_{JX} \Phi) (JY,Z),  \quad \forall X,Y, Z \in \mathfrak{X} (M).
\label{eq:aux1}
\end{equation}
According to Lemma \ref{teor:natural} and  formula (\ref{eq:potential-torsion}), one obtains
\[
(\nabla^{\mathrm{g}}_X \Phi) (Y,Z) =-\frac{1}{2} (g(\mathrm{T}^{\mathrm{w}}(JZ,X),Y)-g(\mathrm{T}^{\mathrm{w}}(Y,X),JZ) + g(\mathrm{T}^{\mathrm{w}}(Z,X),JY)-g(\mathrm{T}^{\mathrm{w}}(JY,X),Z)),
\]
for all vector fields $X, Y, Z$ on $M$.
Condition (\ref{eq:welladapted}) evaluated in $(JZ,X,Y)$ reads as
\[
g(\mathrm{T}^{\mathrm{w}}(JZ,X),Y)-g(\mathrm{T}^{\mathrm{w}}(Y,X),JZ)=g(\mathrm{T}^{\mathrm{w}}(Z,X),JY)-g(\mathrm{T}^{\mathrm{w}}(JY,X),Z),    \quad \forall X,Y, Z \in \mathfrak{X} (M),
\]
and then
\begin{eqnarray*}
(\nabla_X \Phi) (Y,Z) &=&  g(\mathrm{T}^{\mathrm{w}}(JY,X),Z)-g(\mathrm{T}^{\mathrm{w}}(Z,X),JY)),\\
(\nabla_{JX} \Phi) (JY,Z) &=& -\alpha g(\mathrm{T}^{\mathrm{w}}(JY,X),Z)-\alpha g(\mathrm{T}^{\mathrm{w}}(Z,JX),Y),   \quad \forall X,Y, Z \in \mathfrak{X} (M).
\end{eqnarray*}
Taking into account the above equalities and formula (\ref{eq:aux1}) one concludes
\[
\alpha \mathrm{T}^{\mathrm{w}}(Z,X)=J\mathrm{T}^{\mathrm{w}}(Z,JX),   \quad \forall X, Z \in \mathfrak{X} (M).
\]
Then
\[
\mathrm{T}^{\mathrm{w}}(JX,JY)=\alpha\mathrm{T}^{\mathrm{w}}(X,Y),   \quad \forall X,Y \in \mathfrak{X} (M).
\]
The last relation together the hypothesis prove the well adapted connection is torsion-free, i.e.,   $\nabla^{\mathrm{w}}=\nabla^{\mathrm{g}}$, and then $(M,J,g)$ is a K\"{a}hler type manifold. $\blacksquare$
\bigskip

We finish the study of the well adapted connection characterizing  K\"{a}hler type manifolds.

\begin{prop}
Let $(M,J,g)$ be a  $(J^2=\pm1)$-metric manifold. Then $(M,J,g)$ is a K\"{a}hler type manifold if and only if
\[
J\mathrm{T}^{\mathrm{w}}(JX,Y)=\alpha\mathrm{T}^{\mathrm{w}}(X,Y),   \quad \forall X,Y \in \mathfrak{X} (M).
\]
\end{prop}

{\bf Proof}.

$\Rightarrow)$ In this case, $\nabla^{\mathrm{w}}=\nabla^{\mathrm{g}}$, and the condition is trivially satisfied.

$\Leftarrow)$ Taking into account Lemma \ref{teor:primerapropiedad}  and  formula (\ref{eq:welladapted})  reads as
\[
g(\mathrm{T}^{\mathrm{w}}(X,Y),Z)-g(\mathrm{T}^{\mathrm{w}}(Z,Y),X) = -\alpha (g(J\mathrm{T}^{\mathrm{w}}(JX,Y),Z)-g(J\mathrm{T}^{\mathrm{w}}(JZ,Y),X)),   \quad \forall X,Y, Z \in \mathfrak{X} (M),
\]
and by the hypothesis one has
\[
g(\mathrm{T}^{\mathrm{w}}(X,Y),Z)=g(\mathrm{T}^{\mathrm{w}}(Z,Y),X),    \quad \forall X,Y, Z \in \mathfrak{X} (M),
\]
which is equivalent to the vanishing of the torsion tensor of  $\nabla^{\mathrm{w}}$, and then $\nabla^{\mathrm{w}}=\nabla^{\mathrm{g}}$, thus proving the Levi Civita connection of $g$ is a natural connection respect to $(J,g)$. By   Lemma \ref{teor:kahler-caracterizacion}  one concludes  $(M,J,g)$ is a K\"{a}hler type manifold. $\blacksquare$

\subsection{The Kobayashi-Nomizu and the Yano connections}

According to Definitions given by formulas  (\ref{eq:Jconexion1larga}) and (\ref{eq:leyyano2}) one has:

\begin{defin}
Let $(M,J,g)$ be a $(J^2=\pm 1)$-metric manifold.  The Kobayashi-Nomizu connection of $(M,J,g)$ is the linear connection whose covariant derivative is given by
\begin{equation}
\nabla^{\mathrm{kn}}_X Y = \nabla^{0}_ X Y + \frac{(-\alpha)}{4}((\nabla^{\mathrm{g}}_YJ) JX - (\nabla^{\mathrm{g}}_{JY} J) X),   \quad \forall X,Y \in \mathfrak{X} (M).
\label{eq:leykn}
\end{equation}
\end{defin}

\begin{defin}
Let $(M,J,g)$ be a $(J^2=\pm 1)$-metric manifold.  The Yano connection of $(M,J,g)$ is the linear connection whose covariant derivative is given by
\begin{equation}
\nabla^{\mathrm{y}}_X Y = \nabla^{\mathrm{g}}_X Y +\frac{(-\alpha)}{2} (\nabla^{\mathrm{g}}_Y J) JX +  \frac{(-\alpha)}{4}((\nabla^{\mathrm{g}}_X J) JY - (\nabla^{\mathrm{g}}_{JX} J) Y),   \quad \forall X,Y \in \mathfrak{X} (M).
\label{eq:leyyano}
\end{equation}
\end{defin}
 
 \begin{obs} As we know $\nabla^{\mathrm{kn}}$ is always natural respect to $J$. According to Proposition \ref{teor:leyKN-integrable},  $\nabla^{\mathrm{kn}}$ is torsion-free if and only if $J$ is integrable. But, in general, it is not reducible to the $G_{(\alpha,\varepsilon)}$-structure defined by $(J,g)$. Taking into account Lemma \ref{teor:naturalQ},  $\nabla^{\mathrm{kn}}$ is natural respect to $(J,g)$
 if and only if
\begin{equation}
g((\nabla^{\mathrm{g}}_Y J) JX, Z)-g((\nabla^{\mathrm{g}}_{JY} J)X, Z)+ g((\nabla^{\mathrm{g}}_Z J)JX, Y)-g((\nabla^{\mathrm{g}}_{JZ} J) X, Y)=0,   \quad \forall X, Y, Z \in \mathfrak{X} (M).
\label{eq:KN-aeconnection}
\end{equation}

The Yano connection is adapted to the $J$-structure if and only if it is integrable (Proposition \ref{Yanoadaptada}), thus proving it is not  natural respect to $(J,g)$ in general.
\end{obs}

The following result gives a condition about $\nabla^{\mathrm{kn}}=\nabla^{\mathrm{w}}$.

\begin{prop}
\label{teor:nkn=nw}
Let $(M,J,g)$ be a  $(J^2=\pm1)$-metric manifold. If the Kobayashi-Nomizu  connection is reducible to the $G_{(\alpha,\varepsilon)}$-
structure defined by $(J,g)$ then the Kobayashi-Nomizu  connection and the well-adapted connection coincide.
\end{prop}

{\bf Proof.} By Lemmas \ref{teor:JpropiedadKN} and \ref{teor:propiedadesNJ} one obtains the following relation about the torsion tensor of the Kobayashi-Nomizu  connection
\[
 -\varepsilon(g(\mathrm{T}^{\mathrm{kn}}(JX,Y),JZ)-g(\mathrm{T}^{\mathrm{kn}}(JZ,Y),JX)) = g(\mathrm{T}^{\mathrm{kn}}(X,Y),Z)-g(\mathrm{T}^{\mathrm{kn}}(Z,Y),X),   \quad \forall X,Y, Z \in \mathfrak{X} (M).
\]
Then the Kobayashi-Nomizu  connection satisfies the condition (\ref{eq:welladapted}) in Theorem \ref{teor:bienadaptada-ae-estructura}. As the Kobayashi-Nomizu   connection is assumed to be reducible, one also has $\nabla^{\mathrm{kn}} g=0$, and then  $\nabla^{\mathrm{kn}}=\nabla^{\mathrm{w}}$ by Theorem \ref{teor:bienadaptada-ae-estructura}. $\blacksquare$

\begin{corol}
\label{teor:KN-natural}
Let $(M,J,g)$ be a $(J^2=\pm1)$-metric manifold. Then $(M,J,g)$ is a  quasi-K\"{a}hler type manifold if and only if the Kobayashi-Nomizu  connection is natural respect to the $(\alpha,\varepsilon)$-structure $(J,g)$.
\end{corol}

{\bf Proof.}  If $\alpha\varepsilon=1$ then, by relations  (\ref{eq:tensorNJ2}) and  (\ref{eq:tensorNJ4}), formula (\ref{eq:KN-aeconnection}) reads as
\[
-g(\widetilde N_J^{1}(Y,Z),X)=0,    \quad \forall X, Y, Z \in \mathfrak{X} (M),
\]
thus proving the Kobayashi-Nomizu  connection is  metric if and only if $(M,J,g)$ is a quasi-K\"{a}hler  type manifold (according to Proposition \ref{teor:caracterizacion-quasikahler+1}).

In the case $\alpha\varepsilon=-1$ the above quoted formulas (\ref{eq:tensorNJ2}) and (\ref{eq:tensorNJ4}) allow formula  (\ref{eq:KN-aeconnection}) to be read as 
\[
-(g((\nabla^{\mathrm{g}}_ J) JZ, X) -g((\nabla^{\mathrm{g}}_{JY} J)Z, X)  +g((\nabla^{\mathrm{g}}_Z J)JY, X)-g( (\nabla^{\mathrm{g}}_{JZ}J)Y, X)=0, \quad \forall X, Y, Z \in \mathfrak{X} (M),
\]
thus proving the Kobayashi-Nomizu  connection is  metric if and only if $(M,J,g)$ is a quasi-K\"{a}hler  type manifold (according to Remark \ref{teor:ttnj4}). $\blacksquare$
\bigskip

The Kobayashi-Nomizu  connection has been studied in  papers about specific $(J^2=\pm1)$-metric manifolds (see e.g., \cite{chursinetal}, \cite{ganchev-kassabov}, \cite{gribacheva}).
\bigskip

We study now the situation of the Yano connection. As we have said, it is not a natural connection in general. We want to characterize the case $\nabla^{\mathrm{y}}$ is reducible to the $G_{(\alpha,\varepsilon)}$-structure defined by $(J,g)$. Let   $S^{\mathrm{y}}$ be the diference tensor between  $\nabla^{\mathrm{y}}$ and $\nabla^{\mathrm{g}}$, i.e.,
\begin{equation}
\label{eq:syano}
S^{\mathrm{y}} (X,Y) = \frac{(-\alpha)}{2} (\nabla^{\mathrm{g}}_Y J) JX +  \frac{(-\alpha)}{4}((\nabla^{\mathrm{g}}_X J) JY - (\nabla^{\mathrm{g}}_{JX} J) Y), \quad \forall X, Y \in \mathfrak{X}(M).
\end{equation}
Then one has:
\begin{prop}
\label{teor:yanotecnico2}
Let $(M,J,g)$ be a $(J^2=\pm1)$-metric manifold and let $\nabla^{\mathrm{y}}$ be the covariant derivative of the Yano connection. For all vector fields  $X, Y, Z$ on $M$ the following relations hold:
\begin{enumerate}
\renewcommand*{\theenumi}{\roman{enumi})}
\renewcommand*{\labelenumi}{\theenumi}

\item $\nabla^{\mathrm{kn}}_X Y - \nabla^{\mathrm{y}}_X Y = \frac{(-\alpha)}{4} N_J(X,Y)$.

\item $\mathrm{T}^{\mathrm{y}} (X,Y)= \frac{\alpha}{4} N_J(X,Y)$, where $\mathrm{T}^{\mathrm{y}}$ denotes the torsion tensor of $\nabla^{\mathrm{y}}$.

\item $\nabla^{\mathrm{kn}}_X Y= \nabla^{\mathrm{y}}_X Y- \mathrm{T}^{\mathrm{y}} (X,Y)$.

\item $JS^{\mathrm{y}}(X,Y)-S^{\mathrm{y}}(X,JY)-(\nabla^{\mathrm{g}}_X J) Y = \frac{(-\alpha)}{2} N_J(JX,Y)$.

\item $g(S^{\mathrm{y}}(X,Y),Z)+g(S^{\mathrm{y}}(X,Z),Y)=\frac{(-\alpha)}{2} (\alpha  \varepsilon g((\nabla^{\mathrm{g}}_Y J)Z,JX)+g((\nabla^{\mathrm{g}}_Z J)JX,Y))+\frac{\alpha(1+\alpha\varepsilon)}{4} g((\nabla^{\mathrm{g}}_{JX} J)Y,Z)$.

\end{enumerate}

\end{prop}

{\bf Proof.} The first four items are a direct consequence of  Lemma \ref{teor:yanotecnico} applied to the torsion-free covariant derivative $\nabla^{\mathrm{g}}$. In order to prove item $v)$, observe that, according to formula   (\ref{eq:syano}) and Lemma \ref{teor:tensorNJ}, one has
\begin{eqnarray*}
g(S^{\mathrm{y}}(X,Y),Z)+g(S^{\mathrm{y}}(X,Z),Y)&=& \frac{(-\alpha)}{2} g((\nabla^{\mathrm{g}}_Y J) JX, Z) +  \frac{(-\alpha)}{4}( g((\nabla^{\mathrm{g}}_X J) JY,Z) - g((\nabla^{\mathrm{g}}_{JX} J) Y,Z))\\
                                           &+&\frac{(-\alpha)}{2} g((\nabla^{\mathrm{g}}_Z J) JX, Y) +  \frac{(-\alpha)}{4}( g((\nabla^{\mathrm{g}}_X J) JZ,Y) - g((\nabla^{\mathrm{g}}_{JX} J) Z,Y))\\
                                           &=& \frac{(-\alpha)}{2} (\alpha  \varepsilon g((\nabla^{\mathrm{g}}_Y J)Z,JX)+g((\nabla^{\mathrm{g}}_Z J)JX,Y))+\frac{\alpha(1+\alpha\varepsilon)}{4} g((\nabla^{\mathrm{g}}_{JX} J)Y,Z). \ \blacksquare
\end{eqnarray*}

The following results are direct consequences of the above one.

\begin{corol} Let $(M,J,g)$ be a $(J^2=\pm1)$-metric manifold.
\begin{enumerate}
\renewcommand*{\theenumi}{\roman{enumi})}
\renewcommand*{\labelenumi}{\theenumi}

\item The Kobayashi-Nomizu and the Yano connections coincide if and only if the $\alpha$-structure $J$ es integrable.

\item  The Yano connection is torsion-free if and only if  $J$ is integrable.

\item The Yano connection is adapted to the $\alpha$-structure $J$ if and only if $J$ is integrable.

\end{enumerate}
\end{corol}

{\bf Proof.} This is a direct consequence of items $i)$, $ii)$ and $iv)$ of Proposition \ref{teor:yanotecnico2}. $\blacksquare$

\bigskip

\begin{corol} Let $(M,J,g)$ be a $(J^2=\pm1)$-metric manifold. Then the following relations hold:

\begin{enumerate}
\renewcommand*{\theenumi}{\roman{enumi})}
\renewcommand*{\labelenumi}{\theenumi}

\item Assuming $\alpha \varepsilon =1$, the Yano connection is a metric connection if and only if  $(M,J,g)$ is a K\"{a}hler type manifold.

\item Assuming $\alpha \varepsilon =-1$, the Yano connection is a metric connection if and only if  $(M,J,g)$ is a nearly K\"{a}hler type manifold.
\end{enumerate}
\end{corol}

{\bf Proof.}

$i)$ Let us assume $\alpha \varepsilon =1$. Then by item $v)$ of Lemma \ref{teor:yanotecnico2} and by Lemma \ref{teor:tensorpotencial} one knows  that $\nabla^{\mathrm{y}} g=0$ if and only if
\[
g((\nabla^{\mathrm{g}}_{X} J)Y,Z)= g((\nabla^{\mathrm{g}}_Y J)Z,X)+g((\nabla^{\mathrm{g}}_Z J)X,Y), \quad \forall X, Y, Z \in \mathfrak{X} (M),
\]
which is equivalent, according to  Lemma \ref{teor:kahler+1tecnico}, to the condition $(M,J,g)$ is a K\"{a}hler type manifold.

$ii)$ Now let us assume $\alpha \varepsilon =-1$. Then by item  $v)$ of Lemma \ref{teor:yanotecnico2} and by Lemma \ref{teor:tensorpotencial} one knows  that $\nabla^{\mathrm{y}} g=0$ if and only if
\[
g((\nabla^{\mathrm{g}}_Y J)Z,JX)-g((\nabla^{\mathrm{g}}_Z J)JX,Y)=0, \quad \forall X, Y, Z \in \mathfrak{X} (M),
\]
and according to formula (\ref{eq:tensorNJ2})  one has
\[
g((\nabla^{\mathrm{g}}_Y J)Z,JX)+g((\nabla^{\mathrm{g}}_Z J)Y, JX)=0, \quad \forall X, Y, Z \in \mathfrak{X} (M),
\]
thus concluding
\[
(\nabla^{\mathrm{g}}_Y J) Z +(\nabla^{\mathrm{g}}_Z J) Y =0, \quad \forall Y, Z \in \mathfrak{X} (M),
\]
which is an equivalent condition to that of $(M,J,g)$ being a  nearly K\"{a}hler type manifold, because of Lemma \ref{teor:nearlykahler}. $\blacksquare$

\begin{corol}
\label{teor:yano-natural}
 Let $(M,J,g)$ be a $(J^2=\pm1)$-metric manifold. If the Yano connection is reducible to the $(\alpha ,\varepsilon )$-structure $(J,g)$ then $(M,J,g)$ is a  K\"{a}hler type manifold.
\end{corol}

{\bf Proof.} As the Yano connection is reducible one has $\nabla^{\mathrm{y}} J=0$ and $\nabla^{\mathrm{y}} g=0$. As usual, we distinguish two cases:

If $\alpha\varepsilon =1$ then the result follows from the above one. 

Let us assume $\alpha\varepsilon=-1$. Then by condition  $\nabla^{\mathrm{y}} J=0$ one has  $J$ is integrable. Besides, condition  $\nabla^{\mathrm{y}} g =0$ implies  $(M,J,g)$ is a nearly K\"{a}hler type manifold, and by   Corollary \ref{teor:nearly-quasi} it is a quasi- K\"{a}hler type manifold. Finally, according to Proposition \ref{teor:quasikahler-J-integrable}, $(M,J,g)$ is a  K\"{a}hler type manifold. $\blacksquare$
\bigskip

In the almost Norden case, i.e., in the case of manifolds endowed with a $(-1,-1)$-structure, Yano connections are defined in  \cite{teofilova2} and  \cite{teofilova}. As in those papers the structure $J$ is assumed to be integrable, Yano and Kobayashi-Nomizu connections coincide.

\subsection{Connections with totally skew-symmetric torsion}

Connections with totally skew-symmetric torsion have been studied on the different types of $(J^2=\pm1)$-metric manifolds (see, {\em e.g.}, \cite{agricola}, \cite{friedrich}, \cite{ivanov}, \cite{mekerov2}, \cite{mekerov-manev} and  \cite{teofilova}). As in the rest of the paper, we are looking for a unified treatment of the topic. In the present case, our emphasis is focused on the characterization of the existence of a  connection with totally skew-symmetric torsion adapted to the $G_{(\alpha , \varepsilon)}$-structure. We will obtain the following facts:
\begin{enumerate}
\item Such a  characterization.
\item Assuming $\alpha \varepsilon =-1$, a different characterization and the uniqueness of a  natural connection with totally skew-symmetric torsion.
\item Assuming $\alpha \varepsilon =1$, the equivalence between the existence of a natural  connection with totally skew-symmetric torsion and a global property of the manifold: it is  quasi-K\"{a}hler.
\end{enumerate}

First, let us remember:

\begin{defin}
Let $(M,J,g)$ be a $(J^2=\pm1)$-metric manifold. A connection is said to be a connection with totally skew-symmetric torsion if the operator defined as
\[
g(\mathrm{T}^{\mathrm{sk}}(X,Y),Z),   \quad \forall X, Y, Z \in \mathfrak{X} (M),
\]
is a $3$-form on $M$, where $\mathrm{T}^{\mathrm{sk}}$ denotes the torsion tensor of $\nabla^{\mathrm{sk}}$, this being the covariant derivative of the connection.
\end{defin}

As $\mathrm{T}^{\mathrm{sk}}$ is a skew-symmetric tensor, i.e., $\mathrm{T}^{\mathrm{sk}}(X,Y)=-\mathrm{T}^{\mathrm{sk}}(Y,X)$, for all $X, Y$ vector fields on $M$, the above condition is equivalent to the following one
\[
g(\mathrm{T}^{\mathrm{sk}}(X,Y),Z)=-g(\mathrm{T}^{\mathrm{sk}}(Z,Y),X),   \quad \forall X, Y, Z \in \mathfrak{X} (M).
\]
We are interested in natural connections with totally skew-symmetric torsion. In the following result, which follows directly from formula (\ref{eq:potential-torsion}),  we obtain a relationship between the torsion and the potential tensors of such a connection.

\begin{lema}
\label{teor:torsion-sk}
Let $(M,J,g)$ be a $(J^2=\pm1)$-metric manifold and let $\nabla^{\mathrm{sk}}$ be a natural covariant derivative with totally skew-symmetric torsion. Then the  potencial tensor $S^{\mathrm{sk}}$ of $\nabla^{\mathrm{sk}}$ satisfies
\[
S^{\mathrm{sk}}(X,Y)= \frac{1}{2} \mathrm{T}^{\mathrm{sk}}(X, Y),  \quad \forall X, Y \in \mathfrak{X} (M).
\]
\end{lema}

Then, we have: 
\begin{teor}
\label{teor:condicion-sk}
Let $(M,J,g)$ be a $(J^2=\pm1)$-metric manifold. Then the following conditions are equivalent:
\begin{enumerate}
\renewcommand*{\theenumi}{\roman{enumi})}
\renewcommand*{\labelenumi}{\theenumi}

\item There exists a natural covariant derivative $\nabla^{\mathrm{sk}}$ with tollay skew-symmetric torsion.

\item The Nijenhuis tensor is
\begin{equation}
N_J(X,Y)=2 ((\nabla^{\mathrm{g}}_{X} J)JY + (\nabla^{\mathrm{g}}_{JX} J)Y),   \quad \forall X, Y \in \mathfrak{X} (M).
\label{eq:nijenhuis-quasikahler}
\end{equation}

\item 
The following relation holds:
\begin{equation}
(\nabla^{\mathrm{g}}_{X} J)JY + (\nabla^{\mathrm{g}}_{JX} J)Y + (\nabla^{\mathrm{g}}_{Y} J)JX + (\nabla^{\mathrm{g}}_{JY} J)X=0,   \quad \forall X, Y \in \mathfrak{X} (M).
\label{eq:condicion-sk}
\end{equation}
\end{enumerate}
\end{teor}

{\bf Proof.}

$i)\Rightarrow ii)$ According to the above Lemma and Lemma \ref{teor:natural} one has
\[
(J\mathrm{T}^{\mathrm{sk}}(JX, Y) -\alpha \mathrm{T}^{\mathrm{sk}} (X, Y)) + (J\mathrm{T}^{\mathrm{sk}}(X, JY)   - \mathrm{T}^{\mathrm{sk}} (JX, JY)) = 2 ((\nabla^{\mathrm{g}}_{X} J)JY + (\nabla^{\mathrm{g}}_{JX} J)Y),
\]
and by Lemma \ref{teor:nijenhuis-torsion} one obtains the expression of the Nijenhuis tensor.

$ii) \Rightarrow iii)$ Let us consider the expression of the Nijenuis tensor given in $ii)$ and that given in formula (\ref{eq:nijenhuis2}). Subtracting both expressions one has
\[
(\nabla^{\mathrm{g}}_{X} J)JY + (\nabla^{\mathrm{g}}_{JX} J)Y + (\nabla^{\mathrm{g}}_{Y} J)JX + (\nabla^{\mathrm{g}}_{JY} J)X=0,   \quad \forall X, Y \in \mathfrak{X} (M),
\]
as we wanted.

$iii) \Rightarrow i)$ We should define a natural covariant derivative with totally skew-symmetric torsion. Let us consider separately both cases $\alpha\varepsilon=\pm 1$.

In the case $\alpha\varepsilon=-1$ let $\nabla$ be the covariant derivative defined as follows
\begin{equation}
\label{eq:sk-1}
g(\nabla_X Y, Z)= g(\nabla^{0}_X Y, Z) +\frac{(-\alpha)}{2} (g((\nabla^{\mathrm{g}}_Y J) JZ,X)+g((\nabla^{\mathrm{g}}_{JZ} J)Y, X)),   \quad \forall X, Y, Z \in \mathfrak{X} (M).
\end{equation}
As we have defined $\nabla$ from the first canonical connection $\nabla ^{0}$, in order to prove $\nabla$ is natural we must check we are in the conditions of Lemma \ref{teor:naturalQ}. For this, we should consider the tensor $Q=\nabla - \nabla ^{0}$.  We must prove the following two conditions, according to Lemma \ref{teor:naturalQ}: 

\begin{enumerate}
\item $Q\in \mathcal L_{\alpha}$, i.e., $Q(X,JY)=JQ(X,Y)$, 
\item $g(Q(X,Y),Z)+g(Q(X,Z),Y)=0$,
\end{enumerate}
for any vector fields $X, Y \in \mathfrak{X} (M)$. Taking into account the hypothesis $iii)$, one has
\begin{eqnarray*}
g(Q(X,JY)-JQ(X,Y),Z)& =&g(Q(X,JY),Z)-g(JQ(X,Y),Z)\\
&=& \frac{(-\alpha)}{2}  g( \alpha (\nabla^{\mathrm{g}}_{Y} J) Z + (\nabla^{\mathrm{g}}_{JY} J) JZ + (\nabla^{\mathrm{g}}_{JZ} J)JY  + \alpha (\nabla^{\mathrm{g}}_{Z} J)Y, X)=0,\\
g(Q(X,Y),Z)+g(Q(X,Z),Y)&=&  \frac{(-\alpha)}{2}  g( (\nabla^{\mathrm{g}}_{Y} J) JZ + (\nabla^{\mathrm{g}}_{JY} J) Z + (\nabla^{\mathrm{g}}_{Z} J)JY  +  (\nabla^{\mathrm{g}}_{JZ} J)Y, X)=0,
\end{eqnarray*}
for all $X, Y, Z \in \mathfrak{X} (M)$, thus proving $\nabla$   is a natural covariant derivative.

In order to prove the torsion tensor $\mathrm{T}$ of $\nabla$ is totally skew-symmetric, observe that, according to Lemma \ref{teor:tensorNJ} one has
\[
g(\mathrm{T}(X,Y),Z)= (-\alpha) (g((\nabla^{\mathrm{g}}_X J) JY,Z)-g((\nabla^{\mathrm{g}}_Y J) JX,Z)+g((\nabla^{\mathrm{g}}_{JZ} J)Y, X)),   \quad \forall X, Y, Z \in \mathfrak{X} (M),
\]
and then
\[
g(\mathrm{T}(X,Y),Z) + g(\mathrm{T}(Z,Y),X) = \alpha g((\nabla^{\mathrm{g}}_X J) JZ+ (\nabla^{\mathrm{g}}_{JX} J) Z+(\nabla^{\mathrm{g}}_Z J) JX +(\nabla^{\mathrm{g}}_{JZ} J) X,Y)=0, \quad \forall X, Y, Z \in \mathfrak{X} (M),
\]
thus proving $\mathrm{T}$ is totally skew-symmetric.
\bigskip

In the case $\alpha\varepsilon=1$ we use a similar idea, defining  $\nabla$ as the covariant derivative given by
\begin{equation}
\nabla_X Y = \nabla^{0}_X Y + \frac{\alpha}{4} ((\nabla^{\mathrm{g}}_Y J)JX - (\nabla^{\mathrm{g}}_{JY} J)X),   \quad \forall X, Y \in \mathfrak{X} (M).
\label{eq:sk+1}
\end{equation}

In order to prove $\nabla$ is natural we must check the same properties for the  tensor $Q=\nabla - \nabla ^{0}$, according to the same Lemma \ref{teor:naturalQ}. One has
\begin{eqnarray*}
Q(X,JY)- JQ(X,Y)&=& \frac{\alpha}{4} (\nabla^{\mathrm{g}}_{JY} J)JX - \frac{1}{4} (\nabla^{\mathrm{g}}_Y J) X- \left(-\frac{1}{4} (\nabla^{\mathrm{g}}_Y J) X + \frac{\alpha}{4}  (\nabla^{\mathrm{g}}_{JY} J)JX\right)=0,\\
g(Q(X,Y),Z)+g(Q(X,Z),Y)&=& \frac{(-\alpha)}{4} g((\nabla^{\mathrm{g}}_Y J)JZ +(\nabla^{\mathrm{g}}_{JY} J)Z + (\nabla^{\mathrm{g}}_Z J)JY +(\nabla^{\mathrm{g}}_{JZ} J)Y,X)=0,
\end{eqnarray*}
for all $X, Y, Z \in \mathfrak{X} (M)$, thus proving  $\nabla$  is a natural covariant derivative.

About the torsion, from hypothesis $iii)$ one obtains
\[
\mathrm{T}(X,Y)= \frac{3\alpha}{2} (\nabla^{\mathrm{g}}_Y J)JX + \alpha (\nabla^{\mathrm{g}}_{JX} J) Y + \frac{\alpha}{2} (\nabla^{\mathrm{g}}_{JY} J)X,   \quad \forall X, Y \in \mathfrak{X} (M),
\]
and then
\[
g(\mathrm{T}(X,Y),Z) +  g(\mathrm{T}(Z,Y),X) =\alpha(g( (\nabla^{\mathrm{g}}_{JX} J)Y, Z)  + g((\nabla^{\mathrm{g}}_{JY} J) X,Z) + g( (\nabla^{\mathrm{g}}_{JZ} J)X,Y)),    \quad \forall X, Y, Z \in \mathfrak{X} (M).
\]

We must prove the torsion is totally skew-symmetric, i.e.,  $g(\mathrm{T}(X,Y),Z) +  g(\mathrm{T}(Z,Y),X) =0$, for all vector fields on $M$. In this case $\alpha\varepsilon=1$, condition $iii)$ means the second Nijenhuis tensor vanishes, according to formula (\ref{eq:njt}). Then, by   Proposition \ref{teor:caracterizacion-quasikahler+1}  one has $(M,J,g)$ is a  quasi-K\"{a}hler type manifold. Evaluating the condition of a quasi-K\"{a}hler type manifold of Definition \ref{teor:quasikahler+1}  in $(JX, JY, JZ)$, and taking into account formula (\ref{eq:tensorNJ3}), one obtains 
\[
\alpha(g( (\nabla^{\mathrm{g}}_{JX} J)Y, Z)  + g((\nabla^{\mathrm{g}}_{JY} J) X,Z) + g( (\nabla^{\mathrm{g}}_{JZ} J)X,Y))=0,    \quad \forall X, Y, Z \in \mathfrak{X} (M),
\]
thus proving $\mathrm{T}$ is totally skew-symmetric. $\blacksquare$
\bigskip

This theorem provides a common characterization of the existence of a natural connection with totally skew-symmetric torsion for all the four  geometries, which we classify by two conditions:  $\alpha\varepsilon=-1$ and $\alpha\varepsilon=1$. One can recover the specific results for each geometry, as we show in the following Propositions:

\begin{prop}[{\cite[Theor.\ 10.1]{friedrich}, \cite[Corol.\ 3.3]{ivanov}}]
\label{teor:sk-1-nijenhuis}
Let $(M,J,g)$ be a $(J^2=\pm1)$-metric manifold with $\alpha\varepsilon=-1$. There exists a natural connection with totally skew-symmetric torsion if and only if the $(0,3)$ type tensor  defined as
\begin{equation}
g(N_J(X, Y),Z),   \quad \forall X, Y, Z \in \mathfrak{X} (M),
\label{eq:nijenhuis-3forma}
\end{equation}
is a $3$-form on $M$. Besides, this connection is uniquely determined.
\end{prop}

{\bf Proof.} By formula (\ref{eq:nijenhuis2})  and Lemma \ref{teor:tensorNJ} one has
\begin{eqnarray}
g(N_J(X,Y),Z)   &=& - g ((\nabla^{\mathrm{g}}_X J) JZ, Y) - g((\nabla^{\mathrm{g}}_{JX} J) Z,Y) - g((\nabla^{\mathrm{g}}_Y J) JX,Z) - g((\nabla^{\mathrm{g}}_{JY} J) X,Z), \label{eq:nijenhuis-3forma-1}\\
-g(N_J(Z,Y),X)  &=& g ((\nabla^{\mathrm{g}}_Z J) JX, Y) + g((\nabla^{\mathrm{g}}_{JZ} J) X,Y) - g((\nabla^{\mathrm{g}}_Y J) JX,Z) - g((\nabla^{\mathrm{g}}_{JY} J) X,Z). \label{eq:nijenhuis-3forma-2}
\end{eqnarray}
As the Nijenhuis tensor is skew-symmetric (by Lemma \ref{teor:propiedadesNJ}), then property (\ref{eq:nijenhuis-3forma}) is satisfied if and only if
\[
g(N_J(X,Y),Z) = -g(N_J(Z,Y),X),   \quad \forall X, Y, Z \in \mathfrak{X} (M),
\]
which, according to formulas (\ref{eq:nijenhuis-3forma-1}) and (\ref{eq:nijenhuis-3forma-2}),  is satisfied if and only if
\[
(\nabla^{\mathrm{g}}_X J) JZ + (\nabla^{\mathrm{g}}_{JX} J) Z + (\nabla^{\mathrm{g}}_Z J) JX + (\nabla^{\mathrm{g}}_{JZ} J) X =0,   \quad \forall X, Z \in \mathfrak{X} (M), 
\]
which is condition $iii)$ of Theorem \ref{teor:condicion-sk}, and then  this is equivalent to the existence of a natural connection with totally skew-symmetric torsion.

Now, we are going to prove the uniqueness of such a connection. Let  $\nabla^{\mathrm{sk}}$ be an adapted covariant derivative with totally skew-symmetric torsion, and let 
 $\mathrm{T}^{\mathrm{sk}}$ (resp.  $S^{\mathrm{sk}}$) denote its torsion tensor (resp. potential tensor). By Lemma \ref{teor:torsion-sk} one has
\[
S^{\mathrm{sk}}(X,Y)=\frac{1}{2} \mathrm{T}^{\mathrm{sk}}(X,Y), \quad \forall X, Y \in \mathfrak{X} (M).
\]
Besides, according to Lemmas \ref{teor:natural} and \ref{teor:primerapropiedad}, as $\alpha\varepsilon=-1$, one has
\[
g((\nabla^{\mathrm{g}}_X J) Y, Z)= \frac{1}{2} g(J\mathrm{T}^{\mathrm{sk}}(X,Y),Z)-\frac{1}{2} g(\mathrm{T}^{\mathrm{sk}}(X,JY),Z)=-\frac{1}{2} (g(\mathrm{T}^{\mathrm{sk}}(X,Y),JZ)+g(\mathrm{T}^{\mathrm{sk}}(X,JY),Z))
\]
for all vector fields $X,Y, Z$ on $M$. Then
\begin{eqnarray*}
\frac{(-\alpha)}{2} g((\nabla^{\mathrm{g}}_X J)JY, Z)&=& \frac{\alpha}{4} (g(\mathrm{T}^{\mathrm{sk}}(X, JY),JZ)+\alpha g(\mathrm{T}^{\mathrm{sk}}(X,Y),Z)),\\
\frac{(-\alpha)}{2} g((\nabla^{\mathrm{g}}_Y J)JZ, X)&=&\frac{\alpha}{4} (g(\mathrm{T}^{\mathrm{sk}}(Y, JZ),JX)+\alpha g(\mathrm{T}^{\mathrm{sk}}(Y,Z),X)),\\
\frac{(-\alpha)}{2}g((\nabla^{\mathrm{g}}_{JZ} J)Y, X)&=&\frac{\alpha}{4} (g(\mathrm{T}^{\mathrm{sk}}(JZ, Y),JX)+g(\mathrm{T}^{\mathrm{sk}}(JZ,JY),X)), \quad \forall X, Y, Z \in \mathfrak{X} (M).
\end{eqnarray*}
As the torsion tensor of $\nabla^{\mathrm{sk}}$ is skew-symmetric and totally skew-symmetric, the following relations hold:
\begin{eqnarray*}
g(\mathrm{T}^{\mathrm{sk}}(X, JY),JZ)+g(\mathrm{T}^{\mathrm{sk}}(JZ,JY),X)&=&0, \\
g(\mathrm{T}^{\mathrm{sk}}(Y, JZ),JX)+g(\mathrm{T}^{\mathrm{sk}}(JZ, Y),JX)&=&0, \\
 g(\mathrm{T}^{\mathrm{sk}}(X,Y),Z)+g(\mathrm{T}^{\mathrm{sk}}(Y,Z),X)&=&2g(\mathrm{T}^{\mathrm{sk}}(X,Y),Z),  \quad \forall X, Y, Z \in \mathfrak{X} (M).
\end{eqnarray*}
Then one has
\[
\frac{1}{2} g(\mathrm{T}^{\mathrm{sk}}(X,Y),Z)=\frac{(-\alpha)}{2} g((\nabla^{\mathrm{g}}_X J)JY, Z)+ \frac{(-\alpha)}{2} g((\nabla^{\mathrm{g}}_Y J)JZ, X)+ \frac{(-\alpha)}{2}g((\nabla^{\mathrm{g}}_{JZ} J)Y, X),  \quad \forall X, Y, Z \in \mathfrak{X} (M),
\]
thus obtaining the following relation
\[
g(S^{\mathrm{sk}}(X,Y),Z)=\frac{(-\alpha)}{2} g((\nabla^{\mathrm{g}}_X J)JY, Z)+ \frac{(-\alpha)}{2} (g((\nabla^{\mathrm{g}}_Y J)JZ, X)+g((\nabla^{\mathrm{g}}_{JZ} J)Y, X)), \quad \forall X, Y, Z \in \mathfrak{X} (M).
\]
As the above relation is true for all $Z\in \mathfrak{X} (M)$, then the potential tensor $S^{\mathrm{sk}}(X,Y)$ is uniquely determined, and then, the covariant derivative $\nabla^{\mathrm{sk}} =\nabla^{\mathrm{g}} +S^{\mathrm{sk}}$ is also uniquely determined, as we wanted. 

The expression of this unique covariant derivative is that given in  (\ref{eq:sk-1}). $\blacksquare$

\begin{obs} According to Proposition \ref{teor:nearly-quasi} and Theorem \ref{teor:condicion-sk}, one can easily deduce that a nearly K\"{a}hler type manifold has a natural connection with totally skew-symmetric torsion. And, by
 (\ref{eq:sk-1}), one obtains that this unique connection is the first canonical connection
\[
g(\nabla_X Y, Z)= g(\nabla^{0}_X Y, Z) +\frac{(-\alpha)}{2} (g((\nabla^{\mathrm{g}}_Y J) JZ,X)+g((\nabla^{\mathrm{g}}_{JZ} J)Y, X))= g(\nabla^{0}_X Y, Z),   \quad \forall X, Y, Z \in \mathfrak{X} (M),
\]
taking into account property $ii)$ of Lemma \ref{teor:nearlykahler}. One can see this result in, {\em e.g.}, \cite{ivanov} and \cite{vezzoni}.
\end{obs}

In the case $\alpha\varepsilon=1$ one has:

\begin{prop}[{\cite[Theor.\ 3.1]{mekerov2}, \cite[Theor.\ 5, Theor.\ 6]{mekerov-manev}, \cite[Corol.\ 2.4]{teofilova}}]
\label{teor:sk+1-quasikahler}
Let $(M, J, g)$ be a $(J^2=\pm1)$-metric manifold with $\alpha\varepsilon=1$. There exists a natural connection with totally skew-symmetric torsion if and only if $(M,J,g)$ is a quasi-K\"{a}hler type manifold.
\end{prop}

{\bf Proof.}  As we have said in the proof $iii) \Rightarrow i)$ of Theorem \ref{teor:condicion-sk}, in the case $\alpha\varepsilon=1$, condition (\ref{eq:condicion-sk}) means the second Nijenhuis tensor vanishes, i.e., $(M,J,g)$ is a quasi-K\"{a}hler type manifold. $\blacksquare$
\bigskip

The following table summarizes some of the main properties of the distinguished connections on a $(J^{2}=\pm 1)$-metric manifold $(M,J,g)$ studied through the present paper.

\begin{table}[htb]
\begin{center}
\begin{tabular}{|c|c|c|}
 \hline
 \hline
Connection & $\alpha \varepsilon=-1$ & $\alpha \varepsilon=1$ \\
  \hline
  \hline
     $\nabla ^{\mathrm{g}}$ Levi-Civita &  \multicolumn{2}{|c|}{adapted if and only if $(M,J,g)$ is K\"{a}hler type} \\
 \hline
  $\nabla ^{0}$ first canonical  & \multicolumn{2}{|c|}{always adapted} \\
  \hline
  & always adapted and & \\
  $\nabla ^{\mathrm{c}}$ Chern & $\nabla ^{0}=\nabla ^{\mathrm{c}}$ if and only if& there is no\\
    &$(M,J,g)$ is quasi-K\"{a}hler type &  \\
  \hline
    & always adapted and & always adapted and \\
   $\nabla ^{\mathrm{w}}$ well-adapted &  $\nabla ^{0}=\nabla ^{\mathrm{w}}$ if and only if  & $\nabla ^{0}=\nabla ^{\mathrm{w}}$ if and only if\\
     & $(M,J,g)$ is quasi-K\"{a}hler type &  $  J$ is integrable\\
  \hline
   $\nabla ^{\mathrm{kn}}$ Kobayashi-Nomizu &\multicolumn{2}{|c|} {adapted if and only if $(M,J,g)$ is quasi-K\"{a}hler type}  \\
  
  \hline
   $\nabla ^{\mathrm{y}}$ Yano &    \multicolumn{2}{|c|}{adapted if and only if $(M,J,g)$ is K\"{a}hler type} \\
  \hline
   $\nabla ^{\mathrm{sk}}$ totally skew- &   \multicolumn{2}{|c|}{adapted if and only if $(\nabla^{\mathrm{g}}_{X} J)JY+(\nabla^{\mathrm{g}}_{JX} J)Y+(\nabla^{\mathrm{g}}_{Y} J)JX+(\nabla^{\mathrm{g}}_{JY}J)X=0$} \\ 
 symmetric torsion & \multicolumn{2}{|c|}{}\\
  \hline
  \hline
\end{tabular}
\end{center}
\caption{Covariant derivatives on a $(J^2=\pm1)$-metric manifold studied in this section}
\label{table:connections}
\end{table}

\section{The one-parameter family of canonical connections}

Canonical connections on a $(J^{2}=\pm1)$-metric manifold $(M,J,g)$ are a class of significative connections adapted to such a structure.  They were introduced in the two geometries of the case $\alpha \varepsilon =-1$, in \cite[Defin. 2]{gauduchon} for the almost Hermitian geometry and in \cite[Defin. 3.4]{ivanov} for the almost para-Hermitian one. This class of connections consists on a one-parameter family of adapted connections which depends on the first canonical connection $\nabla ^{0}$ and the differential of the fundamental form  $\Phi$ (see \cite[Formula (2.5.4)]{gauduchon} and \cite[Formula (3.13)]{ivanov}).  In the case $(\alpha, \varepsilon)=(-1,1)$, in \cite[Formula (11)]{davidov} it is shown another expression of this family of canonical connections in terms of  $\nabla^{\mathrm{g}}$ and $\nabla^{\mathrm{g}} J$.

A key point is that this family of canonical connections is the affine line determined by the first canonical connection and the Chern connection, i.e., the set 
\begin{equation}
\nabla^{\mathrm{t}}= (1-t) \nabla^{0} + t \nabla^{\mathrm{c}},   \quad \forall t \in \mathbb{R}.
\label{eq:conexionescanonicas}
\end{equation}

As the Chern connection can not be defined in the $\alpha \varepsilon =1$ context, it does not seem possible to define canonical connections in this case. But we will show a way to do it. The idea is the following: in the case $\alpha \varepsilon =-1$, the well adapted connection is also a canonical connection, i.e., is a connection in the line defined in (\ref{eq:conexionescanonicas}). Then this line can be parametrized as 
\[
\nabla^{\mathrm{s}}= (1-s) \nabla^{0} + s \nabla^{\mathrm{w}},   \quad \forall s \in \mathbb{R}.
\]
As the first canonical connection and the well adapted connection can be also defined in the case $\alpha \varepsilon =1$, we are able to define canonical connections on any $(J^{2}=\pm1)$-metric manifold $(M,J,g)$.

In the following result we obtain the one-parameter family of canonical connections. First of all, let us remember the expression of the tensor
\[
{\mathcal F} (\nabla^{0}, X, Y, Z) = g(\mathrm{T}^{0}(X,Y),Z)- g(\mathrm{T}^{0}(Z,Y),X) + \varepsilon (g(\mathrm{T}^{0}(JX,Y),JZ)- g(\mathrm{T}^{0}(JZ,Y),JX)),
\]
for all vector fields $X, Y, Z$ on $M$, given in the proof of Theorem \ref{teor:nabla0=nablaw}. Then we have:

\begin{teor}
\label{teor:familiauniparametrica-1}
Let  $(M,J,g)$ be a  $(J^2=\pm 1)$-metric manifold with  $\alpha\varepsilon=-1$. 

$i)$ The one-parameter family of covariant derivatives on $(M,J,g)$  given by
\begin{equation}
g(\nabla^{\mathrm{s}}_X Y,Z) =  g(\nabla^{0}_X Y, Z)+  \frac{\alpha s}{6} {\mathcal F} (\nabla^{0},Y,X,Z)
                     =g(\nabla^{0}_X Y, Z) +\frac{\alpha s}{12} g(\widetilde N_J^{-1} (Y,Z),X),  \quad \forall X, Y, Z \in \mathfrak{X} (M), \forall s \in \mathbb{R},
\label{eq:familiauniparametrica-1}
\end{equation}
is a family of natural covariant derivatives which contains the covariant derivatives of the first canonical, the Chern and the well adapted connections.

$ii)$ The above family and that of canonical connections coincide. Besides, the family can be parametrized as 
\[
\nabla^{\mathrm{s}}= (1-s) \nabla^{0} + s \nabla^{\mathrm{w}},   \quad \forall s \in \mathbb{R}.
\]
\end{teor}

{\bf Proof.} $i)$ The second equality is direct by formula (\ref{eq:functorial-primeracanonica-1}).  By Lemma \ref{teor:propiedadesNJ0} one obtains
\begin{eqnarray*}
g(Q(X,JY),Z)-g(JQ(X,Y),Z)&=& \frac{\alpha s}{12} g(\widetilde N_J^{-1}(JY,Z),X) + \frac{\alpha s}{12} g(\widetilde N_J^{-1}(Y,JZ),X)=0, \\
g(Q(X,Y),Z)+ g(Q(X,Z),Y)    &=& \frac{\alpha s}{12} g(\widetilde N_J^{-1}(Y,Z),X) + \frac{\alpha s}{12} g(\widetilde N_J^{-1}(Z,Y),X)=0, \quad \forall X, Y, Z \in \mathfrak{X} (M),
\end{eqnarray*}
and then, according to Lemma \ref{teor:naturalQ}, each covariant derivative $\nabla^{\mathrm{s}}$, $s\in \mathbb{R}$, is  adapted to  $(J,g)$.

Obviously $\nabla ^{0}$ belongs to the family, taking $s=0$. In order to prove the covariant derivatives of the Chern and  well adapted connections also belongs to the family we must find other values of $s$ which determine these covariant derivatives. Let  $s\in \mathbb{R}$; then the torsion tensor $\mathrm{T}^{\mathrm{s}}$ of $\nabla^{\mathrm{s}}$ satisfies
\begin{equation}
g(\mathrm{T}^{\mathrm{s}}(X,Y),Z)=g(\mathrm{T}^{0}(X,Y),Z)+\frac{\alpha s}{12} (g(\widetilde N_J^{-1}(Y,Z),X)-g(\widetilde N_J^{-1}(X,Z),Y)),
   \quad \forall X, Y, Z \in \mathfrak{X} (M).  \label{eq:tors-1}
\end{equation}
 For all vector fields $X, Y, Z$ on $M$, according to Lemma \ref{teor:propiedadesNJ0} one has:
 \[
 g(\widetilde N_J^{-1}(Y,Z),X) + \varepsilon g(\widetilde N_J^{-1}(Y,JZ)JX) -(g(\widetilde N_J^{-1}(Y,X),Z)+ \varepsilon g(\widetilde N_J^{-1}(Y,JX),JZ)) = -2g(\widetilde N_J^{-1}(X,Z),Y),
\]
then one obtains
\[
g(\mathrm{T}^{\mathrm{s}}(X,Y),Z)-g(\mathrm{T}^{\mathrm{s}}(Z,Y),X) +\varepsilon (g(\mathrm{T}^{\mathrm{s}}(JX,Y),JZ)-g(\mathrm{T}^{\mathrm{s}}(JZ,Y),JX))= \frac{\alpha}{2}(1-s) g(\widetilde N_J^{-1}(X,Z),Y).
\]
Taking $s=1$, the above expression reads as formula (\ref{eq:welladapted}), thus proving,  by  Proposition \ref{teor:bienadaptada-ae-estructura}, that the corresponding covariant derivative is that of the well adapted connection. 

Now, we are going to prove that the Chern connection corresponds to $s=3$.  For all vector fields $X, Y, Z$ on $M$, according to Lemma \ref{teor:tensorNJ} one obtains
\[
g(\widetilde N_J^{-1}(Y,Z),X)-g(\widetilde N_J^{-1}(X,Z),Y)= g(\widetilde N_J^{-1}(X,Y),Z)+2 (g((\nabla^{\mathrm{g}}_Z J)JX,Y )-g((\nabla^{\mathrm{g}}_{JZ}J) X,Y)),
\]
which together formula  (\ref{eq:tors-1}), formula (\ref{eq:toro-alpha}) in the case $\alpha\varepsilon=-1$, and Lemma \ref{teor:propiedadesNJ0}, leads to the following two expressions
\begin{eqnarray}
g(\mathrm{T}^{\mathrm{s}}(X,Y),Z)&=&g(\mathrm{T}^{0}(X,Y),Z)+\frac{\alpha s}{12} g(\widetilde N_J^{-1}(X,Y),Z)+ \frac{\alpha s}{6}( g((\nabla^{\mathrm{g}}_Z J)JX,Y )-g((\nabla^{\mathrm{g}}_{JZ}J) X,Y)), \label{eq:tors-12}\\
g(\mathrm{T}^{\mathrm{s}}(JX,JY),Z)&=& \alpha g(\mathrm{T}^{0}(X,Y),Z)+\frac{1}{12}(6-s) g(\widetilde N_J^{-1}(X,Y),Z)+ \frac{s}{6}( g((\nabla^{\mathrm{g}}_Z J)JX,Y )-g((\nabla^{\mathrm{g}}_{JZ}J) X,Y)). \nonumber
\end{eqnarray}
According to  Theorem \ref{teor:chern-connection}, the covariant derivative $\nabla^{\mathrm{s}}$ is that of  Chern connection if and only if $s=3$. 
\bigskip

$ii)$ Let  $X, Y, Z$ be vector fields on $M$ y $s \in  \mathbb{R}$. Then, one has
\[
g(\nabla^{\mathrm{s}}_X Y, Z)= (1-s) g(\nabla^{0}_X Y, Z)+ s (g(\nabla^{0}_X Y, Z)+\frac{\alpha}{6}{\mathcal F} (\nabla^{0},Y,X,Z))=(1-s) g(\nabla^{0}_X Y, Z) + s g(\nabla^1_X Y, Z),
\]
thus proving
\[
\nabla^{\mathrm{s}}_X Y  = (1-s) \nabla^{0}_X Y +s \nabla^{\mathrm{w}}_X Y,   \quad \forall X, Y \in \mathfrak{X} (M), \forall s \in \mathbb{R},
\]
and then,
\[
\nabla^{\mathrm{s}}= (1-s) \nabla^{0} + s \nabla^{\mathrm{w}},   \quad \forall s \in \mathbb{R}.
\]
Then one has: 
\begin{itemize}
\item the family of the canonical connections is the affine line $\{ \nabla^{\mathrm{t}}= (1-t) \nabla^{0} + t \nabla^{\mathrm{c}} \colon t\in \mathbb{R}\}$ determined by the first canonical and the Chern connections,
\item the first canonical, the Chern and the well adapted connection belong to the family given in the present theorem,
\item the family of the present theorem is a line of adapted connections, $\{ \nabla^{\mathrm{s}}= (1-s) \nabla^{0} + s \nabla^{\mathrm{w}} \colon   s \in \mathbb{R}\}$,
\end{itemize}
and thus one can conclude that both families of canonical connections and that of the present theorem coincide. $\blacksquare$

\begin{obs} 
(1) The family of canonical connections has been parametrized in two different forms. The relation between parameters is $s=3t$. So, the Chern connection corresponds to $t=1$, as it is expressed in \cite{gauduchon},  and to $s=3$. 

(2) In the case of a manifold endowed with an $\alpha$-structure $J$ we had obtained a family of adapted connections respect to that structure in Proposition \ref{teor:uniparametrica-alpha}. That family does not have to do with the present one. In Example \ref{todaslasfamilias} we will be more accurate about this point.

(3) Taking into account the expression of the second Nijenhuis tensor of $(J,g)$, one obtains the following expression of the Chern connection
\[
 g(\nabla^{\mathrm{c}}_X Y , Z)= g(\nabla^{0}_X Y, Z) +\frac{1}{4} g((\nabla^{\mathrm{g}}_Z J) JY-(\nabla^{\mathrm{g}}_Y J)JZ-(\nabla^{\mathrm{g}}_{JZ} J)Y+(\nabla^{\mathrm{g}}_{JY} J)Z,X),  \quad \forall X, Y, Z \in \mathfrak{X} (M).
\]
In the case $(\alpha,\varepsilon)=(-1,1)$, this is the expression of the Chern connection given in  \cite[Formula (10)]{davidov}.
\end{obs}

\begin{ej}
In \cite{gauduchon} the author points out two other distinguished canonical connections, which corresponds to the values $t=-1$ (or $s=-3$ in our notation) and  $t=\frac{1}{3}$ (or $s=1$). This second one is, as we already know, the well adapted connection $\nabla^{\rm{w}}$. The first one is the so-called Bismut connection $\nabla^{\rm{b}}$, firstly introduced in \cite{bismut},  which can be characterized by the following fact: the tensor 
\[
B(X,Y)=\mathrm{T}^{\mathrm{b}}(X,Y)+\frac{\alpha}{4} N_J(X,Y),   \quad \forall X, Y \in \mathfrak{X} (M),
\]
is totally skew-symmetric, where $\mathrm{T}^{\mathrm{b}}$ is the torsion tensor of $\nabla^{\mathrm{b}}$. Indeed,  let $X, Y, Z$ be vector fields on $M$, by formula (\ref{eq:tors-12}) and  Lemma \ref{teor:tensorNJ} one has
\begin{eqnarray*}
g(\mathrm{T}^{-3}(X,Y),Z)+\frac{\alpha}{4} g(N_J(X,Y),Z)&=& \frac{(-\alpha)}{2}( g(\widetilde N_J^{-1}(X,Y),Z) + g((\nabla^{\mathrm{g}}_Z J)JX,Y )-g((\nabla^{\mathrm{g}}_{JZ}J) X,Y)),\\
g(\mathrm{T}^{-3}(Z,Y),X)+\frac{\alpha}{4} g(N_J(Z,Y),X)&=&  \frac{\alpha}{2}( g(\widetilde N_J^{-1}(X,Y),Z) + g((\nabla^{\mathrm{g}}_Z J)JX,Y )-g((\nabla^{\mathrm{g}}_{JZ}J) X,Y)).
\end{eqnarray*}

Under the assumptions of  Theorem \ref{teor:condicion-sk}, formula (\ref{eq:condicion-sk}) implies the second Nijenhuis tensor of $(M,J,g)$ can be written as
\[
\widetilde N_J^{-1}(Y,Z) = 2( (\nabla^{\mathrm{g}}_Y J) JZ+(\nabla^{\mathrm{g}}_{JZ} J)Y), \quad \forall Y, Z \in \mathfrak{X} (M).
\]
Thus, if $s=-3$ one obtains the connection with totally skew-symmetric connection $\nabla^{\mathrm{sk}}$ introduced in (\ref{eq:sk-1}),
\[
g(\nabla^{-3}_X Y, Z)=g(\nabla^{\mathrm{sk}}_X Y, Z) =g(\nabla^{0}_X Y, Z) +\frac{(-\alpha)}{2} (g((\nabla^{\mathrm{g}}_Y J) JZ,X)+g((\nabla^{\mathrm{g}}_{JZ} J)Y, X)),   \quad \forall X, Y, Z \in \mathfrak{X} (M).
\]
If there exists a connection satisfying the previous conditions then such connection belongs to the family of canonical connections.
 
Also one can conclude that the unique natural connection with totally skew-symmetric torsion in the case $\alpha \varepsilon =-1$, if there exists, is the Bismut connection corresponding to the case $s=-3$ (see Proposition \ref{teor:sk-1-nijenhuis}). For this reason in many papers it is called the natural connection with totally skew-symmetric torsion.
\end{ej}

Taking in mind the case $\alpha \varepsilon =-1$, we can go ahead with that of $\alpha \varepsilon =1$:

\begin{defin}
Let $(M,J,g)$ be a $(J^2=\pm1)$-metric manifold with $\alpha\varepsilon=1$. A canonical connection of $(M,J,g)$ will be any connection reducible to the $G_{(\alpha,\varepsilon)}$-structure defined by $(J,g)$ whose covariant derivative has the form
\[
\nabla^{\mathrm{s}}_X Y = (1-s) \nabla^{0}_X Y + s \nabla^{\mathrm{w}}_X Y,   \quad \forall X, Y \in \mathfrak{X} (M), \forall s \in \mathbb{R}.
\]
\end{defin}

Now, we are looking for an expression of the one-parameter family of canonical connections similar to that obtained in Theorem \ref{teor:familiauniparametrica-1}  $i)$ for the case $\alpha \varepsilon =-1$.

\begin{prop}
\label{teor:familiauniparametrica+1}
Let $(M,J,g)$ be a $(J^2=\pm 1)$-metric manifold with $\alpha\varepsilon=1$. The family of canonical connections of  $(M,J,g)$ is given by the expression 
\begin{equation}
g(\nabla^{\mathrm{s}}_X Y,Z) = g(\nabla^{0}_X Y, Z)+ \frac{\alpha s}{4} {\mathcal F} (\nabla^{0},Y,X ,Z)
                    = g(\nabla^{0}_X Y, Z) +\frac{\alpha s}{8} g(N_J(Y,Z),X),   \quad \forall X, Y, Z \in \mathfrak{X} (M), \forall s \in \mathbb{R}.
\label{eq:familiauniparametrica+1}
\end{equation}
\end{prop}

{\bf Proof.} The second equality is a direct consequence of formula (\ref{eq:functorial-primeracanonica+1}). According to 
 Lemma \ref{teor:propiedadesNJ}, one obtains  
\begin{eqnarray*}
g(Q(X,JY),Z)- g(JQ(X,Y),Z) &=& \frac{\alpha s}{8} g(N_J(JY,Z),X)- \frac{\alpha s}{8} g(N_J(Y,JZ),X) =0,\\
g(Q(X,Y),Z)+ g(Q(X,Z),Y)    &=& \frac{\alpha s}{8} g(N_J(Y,Z),X) + \frac{\alpha s}{8} g(N_J(Z,Y),X)= 0,   \quad \forall X, Y, Z \in \mathfrak{X} (M),
\end{eqnarray*}
thus proving, by Lemma \ref{teor:naturalQ}, each covariant derivative  $\nabla^{\mathrm{s}}$, $s\in \mathbb{R}$, is adapted to  $(J,g)$.

Given $s\in \mathbb{R}$, the torsion of $\nabla^{\mathrm{s}}$  satisfies
\begin{eqnarray}
g(\mathrm{T}^{\mathrm{s}}(X,Y),Z)&=&g(\mathrm{T}^{0}(X,Y),Z)+\frac{\alpha s}{8} (g(N_J(Y,Z),X)-g(N_J(X,Z),Y)), \label{eq:tors+1}\\
g(\mathrm{T}^{\mathrm{s}}(Z,Y),X)&=&g(\mathrm{T}^{0}(Z,Y),X)+\frac{\alpha s}{8} (g(N_J(Y,X),Z)-g(N_J(Z,X),Y)), \label{eq:tors+12}
\end{eqnarray}
for all vector fields $X, Y, Z$ on $M$. Then, by Lemma \ref{teor:propiedadesNJ}, one has:
\[
\alpha g(N_J(X,Y),Z)+\alpha \varepsilon g(N_J(X,JY),JZ) = 0,   \quad \forall X, Y, Z \in \mathfrak{X} (M).
\]
Considering the above equalities, formula (\ref{eq:functorial-primeracanonica+1}) and Lemma \ref{teor:propiedadesNJ} one obtains:
\[
g(\mathrm{T}^{\mathrm{s}}(X,Y),Z)-g(\mathrm{T}^{\mathrm{s}}(Z,Y),X) +\varepsilon (g(\mathrm{T}^{\mathrm{s}}(JX,Y),JZ)-g(\mathrm{T}^{\mathrm{s}}(JZ,Y),JX))= \frac{\alpha}{2}(1- s)  g(N_J(X,Z),Y), \
\]
for all vector fields $X, Y, Z$ on $M$. Taking $s=1$, the above expression reads as formula (\ref{eq:welladapted}), thus proving,  by  Proposition \ref{teor:bienadaptada-ae-estructura}, that the corresponding covariant derivative is that of the well adapted connection. 

The first canonical and the well adapted connections belongs to this one-parameter family of natural connections of $(M,J,g)$, because they correspond to the values $s=0$ and $s=1$. By a similar argument to that used  in the proof of item $ii)$ of Theorem \ref{teor:familiauniparametrica-1}, one can conclude that also in the present case one has
\[
\nabla^{\mathrm{s}}_X Y  = (1-s) \nabla^{0}_X Y +s \nabla^{\mathrm{w}}_X Y,   \quad \forall X, Y \in \mathfrak{X} (M), \forall s \in \mathbb{R}.
\]
thus proving both families, that of the statement and that of the canonical connections, coincide. $\blacksquare$
\bigskip

In the case  $\nabla ^{0}=\nabla ^{\mathrm{w}}$ the line of canonical connections reduces to a point. The more geometric properties the manifold has, the less number of different distinguished connections exist. We finish the paper with three examples of distinguished connections on $(J^{2}=\pm 1)$-metric manifolds with additional properties.

\begin{ej}
Assume $\alpha \varepsilon =1$,  $J$ is  non integrable and $(M,J,g)$ is a quasi-K\"{a}hler type manifold. Then $\nabla ^{0}\neq\nabla ^{\mathrm{w}}$, and the well-adapted and the 
Kobayashi-Nomizu  coincide, by Proposition \ref{teor:nkn=nw} and Corollary \ref{teor:KN-natural}. Besides, there exists a covariant derivative  with totally skew-symmetric torsion (see  Theorem \ref{teor:condicion-sk}). Taking into account the expression (\ref{eq:nijenhuis-quasikahler}) of the Nijenhuis tensor of $J$ one obtains that the family of connections in Proposition \ref{teor:familiauniparametrica+1} can be written as
\[
\nabla^{\mathrm{s}}_X Y = \nabla^{0}_X Y + \frac{(-\alpha) s}{4}  ((\nabla^{\mathrm{g}}_{Y} J)JX - (\nabla^{\mathrm{g}}_{JY} J)X) ,   \quad \forall X, Y \in \mathfrak{X} (M), \forall s \in \mathbb{R}.
\]
Taking $s=-1$ one obtains the connection with totally skew-symmetric connection $\nabla^{\mathrm{sk}}$ introduced in (\ref{eq:sk+1}), 
\[
\nabla^{\mathrm{-1}}_X Y = \nabla^{\mathrm{sk}}_X Y =\nabla^{0}_X Y + \frac{\alpha}{4}  ((\nabla^{\mathrm{g}}_{Y} J)JX - (\nabla^{\mathrm{g}}_{JY} J)X) ,   \quad \forall X, Y \in \mathfrak{X} (M), \forall s \in \mathbb{R}.
\]
thus proving this connection belongs to the family of canonical connections. In fact, it is the unique connection in the family with totally skew-symmetric torsion. We prove this claim. Let  $X, Y, Z \in \mathfrak{X} (M)$. Then
\begin{eqnarray*}
g((\nabla^{\mathrm{g}}_{JY}J)X, Z)&=& -g((\nabla^{\mathrm{g}}_X J)JY, Z) - g((\nabla^{\mathrm{g}}_Z J)JY, X), \\
g(\mathrm{T}^{0}(X,Y),Z)+g(\mathrm{T}^{0}(Z,Y),X)&=& \frac{\alpha}{2} g((\nabla^{\mathrm{g}}_{JY} J)X,Z),\\
g(N_J(Y,Z),X)+g(N_J(Y,X),Z)&=& 4 g((\nabla^{\mathrm{g}}_{JY} J)X,Z).
\end{eqnarray*}
Given $s\in \mathbb{R}$, by (\ref{eq:tors+1}) and (\ref{eq:tors+12}) and taking into account the above identities, the torsion of  $\nabla^{\mathrm{s}}$  satisfies
\[
g(\mathrm{T}^{\mathrm{s}}(X,Y),Z) + g(\mathrm{T}^{\mathrm{s}}(Z,Y),X) = \frac{\alpha}{2} (1+s) g((\nabla^{\mathrm{g}}_{JY} J)X,Z),   \quad \forall X, Y, Z \in \mathfrak{X} (M).
\]
Then  the torsion tensor of $\nabla^{\mathrm{s}}$ is totally skew-symmetric if and only if $s=-1$.
\end{ej}

\begin{ej}
In  \cite[Theor. 3.1]{mekerov2} and  \cite[Theor. 6]{mekerov-manev} the authors study the torsion of natural connections with totally skew-symmetric torsion defined on a $(J^2=\pm1)$-metric manifold with $\alpha \varepsilon =1$. They call such connections as KT-connections and RPT-connections, respectively. They choose a distinguished connection, which is the unique canonical connection with totally skew-symmetric torsion obtained in (\ref{eq:sk+1}).
\end{ej}
\begin{ej}
\label{todaslasfamilias}
Let $(M,J,g)$ be a $(J^{2}=\pm 1)$-metric manifold. Recall that the Kobayashi-Nomizu connection and  the well-adapted connection coincide if and only if $(M,J,g)$ is a quasi-K\"{a}hler type manifold.  In this case the families of canonical connections obtained in Theorem \ref{teor:familiauniparametrica-1} and Proposition \ref{teor:familiauniparametrica+1} are the same that the family of Proposition \ref{teor:uniparametrica-alpha}, introduced in manifolds having an $\alpha$-structure  $J$, assuming the symmetric covariant derivative, choosen in the manifold, is that of the Levi Civita connection of $g$.
\end{ej}

\end{document}